\documentclass[preprint,12pt,sort&compress]{elsarticle}

\usepackage[a4paper,margin=1in]{geometry}
\usepackage[T1]{fontenc}
\usepackage[utf8]{inputenc}
\usepackage{lmodern}
\usepackage{microtype}
\usepackage{amsmath,amssymb,amsfonts,mathtools}
\usepackage{amsthm}
\usepackage{enumitem}
\usepackage{needspace}
\usepackage[colorlinks=true,linkcolor=blue,citecolor=blue,urlcolor=blue]{hyperref}

\numberwithin{equation}{section}

\makeatletter
\renewcommand\section{\@startsection{section}{1}{\z@}%
  {2.4ex \@plus 0.8ex \@minus 0.2ex}%
  {1.2ex \@plus 0.2ex}%
  {\normalfont\Large\bfseries}}
\renewcommand\subsection{\@startsection{subsection}{2}{\z@}%
  {2.0ex \@plus 0.6ex \@minus 0.2ex}%
  {0.9ex \@plus 0.2ex}%
  {\normalfont\large\bfseries}}
\makeatother

\theoremstyle{plain}
\newtheorem{theorem}{Theorem}[section]
\newtheorem{proposition}[theorem]{Proposition}
\newtheorem{lemma}[theorem]{Lemma}

\theoremstyle{definition}
\newtheorem{assumption}[theorem]{Assumption}

\theoremstyle{remark}
\newtheorem{remark}{Remark}[section]

\begin{document}

\begin{frontmatter}

\title{\large Strong Well-Posedness for One-Dimensional Quantile-Dependent SDEs with H\"older-Continuous Diffusion\tnoteref{funding}}
\tnotetext[funding]{This work is supported by the National Natural Science Foundation of China (Grant Nos. 12261038 and 12371152), and the Natural Science Foundation of Jiangxi Province (Grant Nos. 20232BAB201004 and 20242BAB23003).}

\author[aff1]{Jianqing DENG}
\ead{2086095994@qq.com}
\author[aff1]{Hua ZHANG\corref{cor1}}
\cortext[cor1]{Corresponding author.}
\ead{zh860801@163.com}
\affiliation[aff1]{
  organization={School of Statistics and Data Science, Jiangxi University of Finance and Economics},
  city={Nanchang},
  postcode={330013},
  country={P. R. China}
}

\begin{abstract}
\setlength{\parindent}{2em}
\indent We establish strong well-posedness for one-dimensional stochastic differential equations whose drift and diffusion coefficients depend on a prescribed quantile of the solution's marginal law. The diffusion coefficient is uniformly non-degenerate and H\"older continuous in space. The drift is the sum of a non-increasing H\"older-continuous component and a globally Lipschitz component of at most linear growth. Both components are Lipschitz continuous in the quantile variable. The main obstacles are the limited spatial regularity of the diffusion and the quantitative control of quantiles in terms of the $L^1$-distance between densities.

\indent We freeze the quantile path and combine uniform H\"older estimates for the resulting densities with a cancellation argument for the second-order diffusion perturbation. This yields an $L^1$-stability estimate for the frozen density flows and, together with uniform tail bounds and positive lower density bounds, a Volterra-type estimate for the frozen quantile map. An exponentially weighted supremum norm makes this map a strict contraction on $[0,T]$, giving a unique self-consistent quantile path and hence a unique strong solution. When the diffusion coefficient is independent of the quantile variable, the Yamada--Watanabe approximation also yields stability with respect to the initial condition.
\end{abstract}

\begin{keyword}
Quantile-Dependent SDEs,
H\"older-Continuous Diffusion Coefficients,
Strong Well-Posedness,
Stability Estimates,
Density Estimates,
Yamada--Watanabe Approximation
\end{keyword}

\end{frontmatter}

\section{Introduction}

Stochastic differential equations (SDEs, for short) whose coefficients depend on the law of the solution form an important class of nonlinear stochastic systems. They arise naturally in mean-field interaction models, statistical physics, and stochastic control. McKean \cite{McKean1966} initiated the probabilistic study of nonlinear parabolic equations associated with such dynamics. Sznitman \cite{Sznitman1991} and Carmona and Delarue \cite{CarmonaDelarue2018} developed systematic theories for interacting particle systems, propagation of chaos, and mean-field models. Distribution-dependent SDEs with singular or irregular drifts were studied by Huang and Wang \cite{HuangWang2019} and R\"ockner and Zhang \cite{RocknerZhang2021}. Related density-dependent equations were investigated by Wang \cite{Wang2023}. Further work on equations with irregular coefficients includes Bao and Huang \cite{BaoHuang2022}, while Huang and Wang \cite{HuangWang2023} considered equations with H\"older-continuous diffusion coefficients.

Distribution-dependent SDEs describe systems whose coefficients depend on the marginal law of the solution. Density-dependent equations form a subclass in which this dependence acts through the marginal densities \cite{WangRen2024}. Here we study coefficients depending on a prescribed quantile of the time-marginal distribution, as in models driven by a median, a population quantile, or a value-at-risk functional \cite{HuKouritzinZheng2024,Kolokoltsov2013}. Although this involves only one statistical feature of the law, the quantile functional is generally not continuous under weak convergence. We therefore work with densities admitting positive lower bounds on a fixed compact interval containing the relevant quantiles. This allows their difference to be controlled by the $L^1$-distance between the densities.

Quantile-dependent interactions have been studied in nonlinear diffusion models. Kolokoltsov \cite{Kolokoltsov2013} developed a heat-kernel approach based on tail bounds, positive lower density bounds, and $L^1$-stability estimates. Hu, Kouritzin, and Zheng \cite{HuKouritzinZheng2024} established strong existence and pathwise uniqueness for a degenerate finite-dimensional system under a weak H\"ormander condition. F\"ollmer and Schied \cite{FollmerSchied2016} studied quantile and value-at-risk functionals in risk theory.

We adapt Kolokoltsov's quantile fixed-point strategy \cite{Kolokoltsov2013} to one-dimensional equations in which both coefficients may depend on a prescribed quantile. The diffusion is uniformly non-degenerate and $\theta$-H\"older continuous in space for some $\theta\in[1/2,1]$. The drift consists of a non-increasing $\beta$-H\"older continuous component and a globally Lipschitz component of at most linear growth, with both components Lipschitz continuous in the quantile variable.

After freezing the quantile path, we establish uniform tail bounds, positive lower density bounds, and $C_b^\gamma$-estimates for the frozen densities. The density and derivative estimates follow from the flow-centered heat-kernel bounds of Menozzi, Pesce, and Zhang \cite{MenozziPesceZhang2021}. A forward--backward duality argument in the Fokker--Planck framework of Bogachev et al. \cite{BogachevKrylovRocknerShaposhnikov2015}, combined with a cancellation argument for the second-order diffusion perturbation, yields $L^1$-stability of the frozen densities. This gives a Volterra estimate for the frozen quantile map, which becomes a strict contraction under an exponentially weighted supremum norm. Banach's fixed-point theorem then provides the unique self-consistent quantile path. For quantile-independent diffusion, a Yamada--Watanabe approximation yields stability with respect to the initial condition.

Section~\ref{sec:model} presents the model and main results. Section~\ref{sec:proofs} develops the frozen density and quantile estimates, proves strong well-posedness by contraction, and establishes initial-data stability when the diffusion is independent of the quantile variable.

\section{Model, Assumptions, and Main Results}\label{sec:model}

Let $T>0$ be fixed. Let
\[
(\Omega,\mathcal F,(\mathcal F_t)_{t\in[0,T]},\mathbb P)
\]
be a filtered probability space satisfying the usual conditions, and let
$W=(W_t)_{t\in[0,T]}$ be a one-dimensional standard Brownian motion with respect to
$(\mathcal F_t)_{t\in[0,T]}$. Fix a quantile level $\alpha\in(0,1)$.

Let $X=(X_t)_{t\in[0,T]}$ be an adapted process. Whenever the law of $X_t$ is absolutely
continuous with respect to Lebesgue measure, we denote its density by $u_t$. Thus,
\[
\mathcal L(X_t)(dx)=u_t(x)\,dx.
\]
We write $Q_\alpha(u_t)$ for the corresponding $\alpha$-quantile, whose precise definition
is given below. We consider the following one-dimensional quantile-dependent SDE:
\begin{equation}
 dX_t=b\bigl(X_t,Q_\alpha(u_t)\bigr)\,dt
 +\sigma\bigl(X_t,Q_\alpha(u_t)\bigr)\,dW_t,
 \qquad X_0=\xi,\quad t\in[0,T].
\label{eq:2.1}
\end{equation}

For a probability density $u$ on $\mathbb R$, define its distribution function by
\[
F_u(x)=\int_{-\infty}^x u(y)\,dy,\qquad x\in\mathbb R.
\]
Its $\alpha$-quantile is defined by
\begin{equation}
Q_\alpha(u)=\inf\{x\in\mathbb R:F_u(x)\ge \alpha\}.
\label{eq:2.2}
\end{equation}
Since $F_u$ is continuous, the definition in \eqref{eq:2.2} implies that
\[
F_u\bigl(Q_\alpha(u)\bigr)=\alpha.
\]
We use this definition throughout the paper. We consider solutions whose laws admit
densities for every $t\in[0,T]$. The precise class of solutions considered in this paper is
specified below.

For later use, define the class of admissible densities by
\begin{equation}
\mathcal D
:=
\left\{
 u\in C(\mathbb R)\cap L^1(\mathbb R)\cap L^\infty(\mathbb R):
 u(x)>0\ \text{for all }x\in\mathbb R,\quad
 \int_{\mathbb R}u(x)\,dx=1
\right\}.
\label{eq:2.3}
\end{equation}
For every $u\in\mathcal D$ and $x<y$,
\[
F_u(y)-F_u(x)=\int_x^y u(r)\,dr>0,
\]
since $u(r)>0$ for every $r\in\mathbb R$. Hence the identity
\[
F_u\bigl(Q_\alpha(u)\bigr)=\alpha
\]
uniquely determines $Q_\alpha(u)$.

Whenever it exists, the density $u_t$ satisfies the following nonlinear Fokker--Planck
equation in the distributional sense:
\begin{equation}
\partial_tu_t(x)
=
\frac12\partial_{xx}\!\left(\sigma^2\bigl(x,Q_\alpha(u_t)\bigr)u_t(x)\right)
-
\partial_x\!\left(b\bigl(x,Q_\alpha(u_t)\bigr)u_t(x)\right),
\qquad t\in(0,T].
\label{eq:2.4}
\end{equation}
Eq.~\eqref{eq:2.4} is the forward equation associated with Eq.~\eqref{eq:2.1}. Its
frozen version will be used in the forward--backward duality argument in
Section~\ref{sec:proofs}.

We next introduce the frozen equation. Let
\[
C_T=C([0,T];\mathbb R)
\]
be the space of continuous real-valued functions on $[0,T]$, equipped with the norm
\[
\|\omega\|_T=\sup_{0\le t\le T}|\omega_t|.
\]
For a deterministic path $\omega\in C_T$, consider the frozen SDE
\begin{equation}
 dX_t^\omega=b(X_t^\omega,\omega_t)\,dt
 +\sigma(X_t^\omega,\omega_t)\,dW_t,
 \qquad X_0^\omega=\xi.
\label{eq:2.5}
\end{equation}
We denote the density of $X_t^\omega$ by $u_t^\omega$. Its existence is proved in
Proposition~\ref{prop:frozen-density}. In the distributional sense, $u_t^\omega$ satisfies
\begin{equation}
\partial_tu_t^\omega(x)
=
\frac12\partial_{xx}\!\left(\sigma^2(x,\omega_t)u_t^\omega(x)\right)
-
\partial_x\!\left(b(x,\omega_t)u_t^\omega(x)\right).
\label{eq:2.6}
\end{equation}
For $\omega\in C_T$, define the frozen quantile map by
\[
(\Phi(\omega))_t:=Q_\alpha(u_t^\omega),\qquad t\in[0,T].
\]
In Section~\ref{sec:proofs}, we prove that $\Phi$ is well defined as a mapping from $C_T$
into itself.
A fixed point of $\Phi$ then gives the self-consistent quantile path associated with
Eq.~\eqref{eq:2.1}.

To control the quantile map, let $K>0$, $\delta>0$, and
\[
0<\varepsilon<\min\{\alpha,1-\alpha\}.
\]
Define
\begin{equation}
\mathcal S_{\alpha,K,\delta,\varepsilon}
:=
\left\{
 u\in\mathcal D:
 \int_{|x|\ge K}u(x)\,dx\le\varepsilon,\quad
 |Q_\alpha(u)|\le K,\quad
 \inf_{|x|\le K}u(x)\ge\delta
\right\}.
\label{eq:2.7}
\end{equation}
The tail condition and $\varepsilon<\min\{\alpha,1-\alpha\}$ imply $|Q_\alpha(u)|<K$. The explicit quantile bound is retained for ease of reference. The positive lower bound on $[-K,K]$ provides a uniform $L^1$-Lipschitz estimate for the quantile map.

We also use the standard bounded H\"older spaces. For $\rho\in(0,1]$ and bounded
$f:\mathbb R\to\mathbb R$, define
\[
[f]_\rho
:=
\sup_{x\ne y}\frac{|f(x)-f(y)|}{|x-y|^\rho},
\qquad
\|f\|_{C_b^\rho}
:=
\|f\|_\infty+[f]_\rho.
\]
We denote by $C_b^\rho(\mathbb R)$ the space of bounded functions with finite
$C_b^\rho$-norm.

To obtain the above bounds uniformly for the frozen equations, we impose the following
assumptions on the coefficients and the initial condition.

\begin{assumption}[Drift Coefficient]\label{ass:b}
The drift coefficient $b:\mathbb R\times\mathbb R\to\mathbb R$ admits the decomposition
\[
b(x,z)=b_1(x,z)+b_2(x,z),\qquad x,z\in\mathbb R.
\]
There exist constants $K>0$ and $\beta\in(0,1]$ such that the following conditions hold
for all $x,y,z,\widetilde z\in\mathbb R$.

\begin{enumerate}
\item[(i)] For each fixed $z\in\mathbb R$, the map $x\mapsto b_1(x,z)$ is non-increasing:
\[
(x-y)\bigl(b_1(x,z)-b_1(y,z)\bigr)\le0.
\]
Moreover, $b_1$ is $\beta$-H\"older continuous in the spatial variable and Lipschitz continuous in the quantile variable:
\[
|b_1(x,z)-b_1(y,z)|\le K|x-y|^\beta,
\qquad
|b_1(x,z)-b_1(x,\widetilde z)|\le K|z-\widetilde z|.
\]

\item[(ii)] $b_2$ is jointly Lipschitz continuous:
\[
|b_2(x,z)-b_2(y,\widetilde z)|
\le K\bigl(|x-y|+|z-\widetilde z|\bigr).
\]

\item[(iii)] The full drift is uniformly bounded at the spatial origin:
\[
|b(0,z)|\le K,\qquad z\in\mathbb R.
\]
\end{enumerate}
\end{assumption}

\begin{remark}
The monotonicity of $b_1$ is used in the Yamada--Watanabe argument, while Lipschitz continuity in the quantile variable allows comparison of frozen equations. The assumptions also yield a uniform linear-growth bound: for all $x,z\in\mathbb R$,
\begin{align*}
|b(x,z)|
&\le |b(0,z)|
  +|b_1(x,z)-b_1(0,z)|
  +|b_2(x,z)-b_2(0,z)|\\
&\le K+K|x|^\beta+K|x|
\le C(1+|x|),
\end{align*}
since $0<\beta\le1$.
\end{remark}

\begin{assumption}[Diffusion Coefficient]\label{ass:sigma}
Let $\sigma:\mathbb R\times\mathbb R\to\mathbb R$. There exist constants
$0<\lambda\le\Lambda<\infty$, $K>0$, and $\theta\in[1/2,1]$ such that, for all
$x,y,z,\widetilde z\in\mathbb R$,
\[
\lambda\le|\sigma(x,z)|\le\Lambda,
\]
and
\[
|\sigma(x,z)-\sigma(y,z)|\le K|x-y|^\theta.
\]
Fix
\[
\varrho:=\min\{\beta,\theta\},
\qquad
0<\gamma<\varrho.
\]
Assume moreover that the map $z\mapsto\sigma(\cdot,z)$ is Lipschitz
continuous as a map into $C_b^\gamma(\mathbb R)$, namely,
\begin{equation*}
\|\sigma(\cdot,z)-\sigma(\cdot,\widetilde z)\|_{C_b^\gamma}
\le K|z-\widetilde z|,
\qquad z,\widetilde z\in\mathbb R.
\end{equation*}
\end{assumption}

\begin{remark}
The condition $\theta\ge1/2$ ensures pathwise uniqueness for the one-dimensional frozen equations, while non-degeneracy and boundedness give the heat-kernel estimates. Up to a change of constant, the $C_b^\gamma$-valued Lipschitz condition combines a uniform Lipschitz bound in the quantile variable with a $\gamma$-H\"older bound for the corresponding spatial increments. It controls the second-order perturbation in the comparison of frozen densities.
\end{remark}

\begin{assumption}[Initial Condition]\label{ass:init}
Let $p\ge2$. The random variable $\xi$ is $\mathcal F_0$-measurable, satisfies
\[
\mathbb E[|\xi|^p]<\infty,
\]
and its law admits a density
\[
u_0\in\mathcal D\cap C_b^\gamma(\mathbb R).
\]
\end{assumption}

\begin{remark}
Since $u_0\in\mathcal D$ is continuous and strictly positive, $Q_\alpha(u_0)$ is unique and $u_0$ has a positive lower bound on every compact interval. Proposition~\ref{prop:density-class} provides constants $K_0,\delta_0>0$ and $0<\varepsilon_0<\min\{\alpha,1-\alpha\}$, independent of $t$ and $\omega$, such that
\[
u_t^\omega\in\mathcal S_{\alpha,K_0,\delta_0,\varepsilon_0},
\qquad t\in[0,T],\quad \omega\in C_T.
\]
The $C_b^\gamma$-regularity of $u_0$ yields the corresponding uniform estimate for the frozen densities.
\end{remark}

We next specify the class of solutions considered in this paper.

An $(\mathcal F_t)$-adapted continuous process $X=(X_t)_{t\in[0,T]}$ is called a strong
solution of Eq.~\eqref{eq:2.1} if the following conditions hold:
\begin{enumerate}
\item[(i)] for every $t\in[0,T]$, the law $\mathcal L(X_t)$ has a density
$u_t\in\mathcal D$,
\item[(ii)] the path
\[
t\longmapsto Q_\alpha(u_t)
\]
belongs to $C_T$, and
\item[(iii)]
\[
X_t
=
\xi
+
\int_0^t b\bigl(X_s,Q_\alpha(u_s)\bigr)\,ds
+
\int_0^t \sigma\bigl(X_s,Q_\alpha(u_s)\bigr)\,dW_s
\]
holds almost surely for every $t\in[0,T]$.
\end{enumerate}

Pathwise uniqueness is said to hold for Eq.~\eqref{eq:2.1} if, whenever $X$ and
$\widetilde X$ are two strong solutions defined on the same filtered probability space,
driven by the same Brownian motion, and satisfying
\[
X_0=\widetilde X_0=\xi\qquad\text{almost surely},
\]
then
\[
\mathbb P\bigl(X_t=\widetilde X_t\ \text{for all }t\in[0,T]\bigr)=1.
\]
Equivalently, $X$ and $\widetilde X$ are indistinguishable. Throughout the paper, pathwise
uniqueness is understood within the above class of strong solutions. We say that
Eq.~\eqref{eq:2.1} is strongly well posed if it admits a strong solution and pathwise
uniqueness holds.

We now state the main results.

\begin{theorem}[Strong Well-Posedness]\label{thm:wellposed}
Suppose that Assumptions~\ref{ass:b}--\ref{ass:init} hold. Then Eq.~\eqref{eq:2.1} is
strongly well posed.
Let $X=(X_t)_{t\in[0,T]}$ denote its unique strong solution. For every $t\in[0,T]$, the law
of $X_t$ admits a density $u_t\in\mathcal D$, and the associated quantile path
\[
q=(Q_\alpha(u_t))_{t\in[0,T]}
\]
belongs to $C_T$.

Furthermore, the following estimate holds for some constant $C>0$:
\begin{equation}
\mathbb E\left[\sup_{0\le t\le T}|X_t|^p\right]
\le
C\bigl(1+\mathbb E[|\xi|^p]\bigr).
\label{eq:2.8}
\end{equation}
The constant $C$ may depend on $p$, $T$, and the constants appearing in
Assumptions~\ref{ass:b} and~\ref{ass:sigma}.
\end{theorem}

\begin{remark}
Theorem~\ref{thm:wellposed} allows both coefficients to depend on the prescribed quantile. The diffusion need only be $\theta$-H\"older continuous in space for some $\theta\in[1/2,1]$, and the drift may combine a non-increasing H\"older-continuous component with a globally Lipschitz component of at most linear growth.
\end{remark}

We next give an initial-data stability result for the case where the diffusion coefficient is
independent of the quantile variable.

\begin{theorem}[Stability with Respect to the Initial Condition]\label{thm:stability}
Suppose that Assumptions~\ref{ass:b} and~\ref{ass:sigma} hold and, in addition,
\[
\sigma(x,z)=\sigma(x),\qquad x,z\in\mathbb R.
\]
Let $\xi$ and $\eta$ be two $\mathcal F_0$-measurable random variables defined on the same
probability space. Suppose that, for some $p\ge2$,
\[
\mathbb E[|\xi|^p]+\mathbb E[|\eta|^p]<\infty.
\]
Assume further that the laws of $\xi$ and $\eta$ admit densities
\[
u_0^\xi,u_0^\eta\in\mathcal D\cap C_b^\gamma(\mathbb R).
\]
Let $X^\xi$ and $X^\eta$ be the unique strong solutions given by
Theorem~\ref{thm:wellposed}, driven by
the same Brownian motion $W$. Denote their densities by $u_t^\xi$ and $u_t^\eta$,
respectively, and define
\[
q_t^\xi=Q_\alpha(u_t^\xi),
\qquad
q_t^\eta=Q_\alpha(u_t^\eta).
\]
Then the following estimates hold for some constant $C>0$:
\begin{equation}
\sup_{0\le t\le T}\mathbb E[|X_t^\xi-X_t^\eta|]
\le
C\mathbb E[|\xi-\eta|].
\label{eq:2.9}
\end{equation}
Moreover, for every $t\in(0,T]$,
\begin{equation}
\|u_t^\xi-u_t^\eta\|_{L^1}
\le
Ct^{-1/2}\mathbb E[|\xi-\eta|],
\label{eq:2.10}
\end{equation}
and
\begin{equation}
|q_t^\xi-q_t^\eta|
\le
Ct^{-1/2}\mathbb E[|\xi-\eta|].
\label{eq:2.11}
\end{equation}
The constant $C$ may depend on $p$, $T$, $\alpha$, and the constants appearing in
Assumptions~\ref{ass:b} and~\ref{ass:sigma}. It may also depend on the moments
$\mathbb E[|\xi|^p]$, $\mathbb E[|\eta|^p]$ and the initial densities
$u_0^\xi$, $u_0^\eta$.
\end{theorem}

\section{Proofs of the Main Results}\label{sec:proofs}

We prove Theorems~\ref{thm:wellposed} and~\ref{thm:stability} by comparing frozen density flows and then controlling the associated quantile map. The main technical point is an integrable estimate for the second-order diffusion perturbation, which turns the quantile stability estimate into a Volterra inequality. An exponentially weighted supremum norm then yields contraction on the whole interval $[0,T]$.

Throughout this section, $C>0$ denotes a generic constant whose value may change from line
to line. Unless otherwise stated, it may depend on $T$, $\alpha$, and the constants in
Assumptions~\ref{ass:b} and~\ref{ass:sigma}. It may also depend on the relevant initial
density or densities and, for moment estimates, on $p$. The constant is independent of the
frozen paths under consideration. A constant used for a positive lower bound is denoted by
$c$, with a subscript added when necessary.

\subsection{Yamada--Watanabe Approximation}

We use the standard Yamada--Watanabe smooth approximation of the absolute value. See
\cite{IkedaWatanabe1989,YamadaWatanabe1971} for the classical construction and
\cite{HuangWang2023} for a related application to SDEs with
H\"older-continuous diffusion coefficients. Throughout this subsection, we write $\sigma(x)$
for the diffusion coefficient.

Fix $r>1$ and $\varepsilon\in(0,1)$. Since
\[
\int_{\varepsilon/r}^{\varepsilon}\frac{1}{x}\,dx=\log r,
\]
there exists a nonnegative function
\[
\psi_{r,\varepsilon}\in C_c^\infty((\varepsilon/r,\varepsilon))
\]
such that
\[
\int_{\varepsilon/r}^{\varepsilon}\psi_{r,\varepsilon}(x)\,dx=1,
\qquad
0\le\psi_{r,\varepsilon}(x)\le\frac{2}{x\log r}.
\]
For $x\in\mathbb R$, define
\[
V_{r,\varepsilon}(x)
=
\int_0^{|x|}\int_0^y\psi_{r,\varepsilon}(z)\,dz\,dy.
\]
Then $V_{r,\varepsilon}\in C^2(\mathbb R)$, and
\begin{equation}
|x|-\varepsilon\le V_{r,\varepsilon}(x)\le|x|,
\qquad x\in\mathbb R.
\label{eq:3.2}
\end{equation}
Moreover,
\begin{equation}
\operatorname{sgn}(x)V_{r,\varepsilon}'(x)\in[0,1],
\qquad x\ne0.
\label{eq:3.3}
\end{equation}
We set $V_{r,\varepsilon}'(0)=0$. Consequently,
\[
|V_{r,\varepsilon}'(x)|\le1,\qquad x\in\mathbb R.
\]
Furthermore,
\[
0\le V_{r,\varepsilon}''(x)
\le
\frac{2}{|x|\log r}\mathbf 1_{[\varepsilon/r,\varepsilon]}(|x|),
\qquad x\ne0.
\]
We also set $V_{r,\varepsilon}''(0)=0$.

By Assumption~\ref{ass:sigma},
\[
|\sigma(x)-\sigma(y)|\le K|x-y|^\theta,
\qquad \theta\in[1/2,1].
\]
For $x\ne y$, we have
\begin{equation}
\begin{aligned}
V_{r,\varepsilon}''(x-y)|\sigma(x)-\sigma(y)|^2
&\le
\frac{2K^2}{\log r}|x-y|^{2\theta-1}
\mathbf 1_{[\varepsilon/r,\varepsilon]}(|x-y|)\\
&\le
\frac{2K^2}{\log r}\varepsilon^{2\theta-1}.
\end{aligned}
\label{eq:3.5}
\end{equation}
Here the last inequality uses $2\theta-1\ge0$. When $x=y$, the left-hand side is zero. If
$\theta>1/2$, then, for fixed $r>1$, the right-hand side converges to zero as
$\varepsilon\downarrow0$. In the critical case $\theta=1/2$, the bound becomes
$2K^2/\log r$. After sending $\varepsilon\downarrow0$, one must also let $r\to\infty$.
The estimate in \eqref{eq:3.5} will be used to control the diffusion term in It\^o's formula
in the proof of Theorem~\ref{thm:stability}.

\subsection{Frozen Equations and Density Estimates}

For a deterministic path $\omega\in C_T$, recall the frozen equation in
Eq.~\eqref{eq:2.5}. Set
\[
b^\omega(t,x):=b(x,\omega_t),
\qquad
\sigma^\omega(t,x):=\sigma(x,\omega_t).
\]
Since $\omega$ is deterministic, this is a time-inhomogeneous one-dimensional SDE with
prescribed coefficients. For $0\le s\le T$ and $x\in\mathbb R$, let
$X_{s,r}^{\omega,x}$ be the solution of
\begin{equation}
\begin{cases}
dX_{s,r}^{\omega,x}
=b(X_{s,r}^{\omega,x},\omega_r)\,dr
+\sigma(X_{s,r}^{\omega,x},\omega_r)\,dW_r,\\
X_{s,s}^{\omega,x}=x,
\end{cases}
\qquad r\in[s,T].
\label{eq:3.6}
\end{equation}
The associated transition operator is defined by
\[
P_{s,t}^\omega f(x)=\mathbb E[f(X_{s,t}^{\omega,x})],
\qquad 0\le s\le t\le T.
\]
For $c>0$, set
\[
g_c(t,z)
:=
t^{-1/2}\exp\!\left(-\frac{|z|^2}{ct}\right),
\qquad t>0,\ z\in\mathbb R.
\]
Let $\chi_{r,s}^\omega$ denote the regularized flow associated with $b^\omega$, constructed using a fixed unit-scale spatial mollification as in \cite[Section~1.2, equations (1.8) and (1.11)]{MenozziPesceZhang2021}. For the lower Gaussian bound alone, \cite[Remark~1.3]{MenozziPesceZhang2021} also allows any Peano flow of the original drift. When this generality is needed, we retain the same notation and explicitly choose $\chi_{r,s}^\omega$ as a fixed solution of
\begin{equation}
\chi_{r,s}^\omega(x)
=x+\int_s^r b(\chi_{\ell,s}^\omega(x),\omega_\ell)\,d\ell.
\label{eq:peano-flow}
\end{equation}

By the uniform spatial derivative estimate for the regularized drift established in
\cite[(1.9)]{MenozziPesceZhang2021}, Gronwall's inequality applied to the forward
and inverse flows yields the following uniform bi-Lipschitz estimate:
there exists a constant $C\ge1$,
independent of $\omega\in C_T$, such that
\begin{equation}
C^{-1}|x-y|
\le |\chi_{r,s}^\omega(x)-\chi_{r,s}^\omega(y)|
\le C|x-y|,
\quad x,y\in\mathbb R,\quad 0\le s\le r\le T.
\label{eq:flow-bilipschitz}
\end{equation}

We next collect the transition-density estimates for the frozen equation that
will be used throughout the sequel. Under Assumptions~\ref{ass:b} and \ref{ass:sigma},
the frozen coefficients satisfy the conditions required for the Gaussian
heat-kernel estimates in \cite[Theorem~1.2]{MenozziPesceZhang2021}.
Consequently, the transition density $p^\omega(s,x;t,y)$ exists and, for
suitable constants $C,c_1,c_2>0$,
\begin{equation}
\begin{aligned}
&C^{-1}g_{c_1}\bigl(t-s,\chi_{t,s}^\omega(x)-y\bigr)\\
&\qquad\le p^\omega(s,x;t,y)
\le Cg_{c_2}\bigl(t-s,\chi_{t,s}^\omega(x)-y\bigr),
\end{aligned}
\label{eq:frozen-gaussian-two-sided}
\end{equation}
for $0\le s<t\le T$ and $x,y\in\mathbb R$. Furthermore,
\begin{equation}
|\partial_xp^\omega(s,x;t,y)|
\le C(t-s)^{-1/2}
 g_c\bigl(t-s,\chi_{t,s}^\omega(x)-y\bigr),
\label{eq:frozen-gradient-kernel}
\end{equation}
and
\begin{equation}
|\partial_{xx}p^\omega(s,x;t,y)|
\le C(t-s)^{-1}
 g_c\bigl(t-s,\chi_{t,s}^\omega(x)-y\bigr).
\label{eq:frozen-second-derivative-kernel}
\end{equation}
The constants in \eqref{eq:frozen-gaussian-two-sided}--\eqref{eq:frozen-second-derivative-kernel}
can be chosen independently of $\omega\in C_T$.
Moreover, by \cite[Remark~1.3]{MenozziPesceZhang2021}, the lower Gaussian
bound in \eqref{eq:frozen-gaussian-two-sided} remains valid when its center
is taken along any Peano flow.
The estimates for the original coefficients follow by the spatial approximation
and compactness argument in
\cite[Section~3.4, (3.13)--(3.14), and Appendix~A]{MenozziPesceZhang2021}.
The corresponding estimates for the mollified coefficients are uniform with
respect to the approximation parameter.

The following proposition collects the main properties and estimates for the frozen equation
that will be used below. In particular, all constants in these estimates can be chosen
independently of the frozen path $\omega$.

\begin{proposition}[Frozen Equation and Density Estimates]\label{prop:frozen-density}\leavevmode\par\noindent
Suppose that Assumptions~\ref{ass:b}--\ref{ass:init} hold. Then, for every
$\omega\in C_T$, the following
assertions hold.

\begin{enumerate}
\item[(i)] The frozen equation in Eq.~\eqref{eq:2.5} admits a unique strong solution
$X^\omega$, and
\begin{equation}
\sup_{\omega\in C_T}
\mathbb E\left[\sup_{0\le t\le T}|X_t^\omega|^p\right]
\le
C\bigl(1+\mathbb E[|\xi|^p]\bigr).
\label{eq:3.7}
\end{equation}

\item[(ii)] The transition process in \eqref{eq:3.6} has a continuous density $p^\omega$ satisfying
\[
P_{s,t}^\omega f(x)
=\int_{\mathbb R}f(y)p^\omega(s,x;t,y)\,dy.
\]
For every bounded measurable function $f$, the transition operator satisfies
\begin{equation}
\|\partial_xP_{s,t}^\omega f\|_{L^\infty}
\le
\frac{C}{\sqrt{t-s}}\|f\|_{L^\infty},
\qquad 0\le s<t\le T.
\label{eq:3.8}
\end{equation}
Moreover, there are constants $C,c>0$, independent of $\omega$, such that
\[
|\partial_{xx}p^\omega(s,x;t,y)|
\le C(t-s)^{-1}
g_c\bigl(t-s,\chi_{t,s}^\omega(x)-y\bigr),
\qquad 0\le s<t\le T,\quad x,y\in\mathbb R.
\]
The constant in \eqref{eq:3.8} is also independent of $\omega$.

\item[(iii)] For every $t\in(0,T]$, the law of $X_t^\omega$ admits a continuous and strictly
positive density $u_t^\omega$, with $u_0^\omega=u_0$. Moreover,
\[
u^\omega\in C([0,T];L^1(\mathbb R)),
\]
and
\begin{equation}
\sup_{\omega\in C_T}\sup_{0\le t\le T}
\|u_t^\omega\|_{C_b^\gamma}
\le C.
\label{eq:3.10}
\end{equation}
\end{enumerate}
\end{proposition}

\begin{proof}
\textup{(i)} By Assumption~\ref{ass:b}, uniformly in $t\in[0,T]$ and $\omega\in C_T$,
\[
|b^\omega(t,x)|\le C(1+|x|).
\]
Moreover, the monotonicity of $b_1$ and the spatial Lipschitz continuity of $b_2$ imply
\[
\begin{aligned}
&(x-y)\bigl(b^\omega(t,x)-b^\omega(t,y)\bigr)\\
&\qquad=(x-y)\bigl(b_1(x,\omega_t)-b_1(y,\omega_t)\bigr)
 +(x-y)\bigl(b_2(x,\omega_t)-b_2(y,\omega_t)\bigr)\\
&\qquad\le K|x-y|^2.
\end{aligned}
\]
On the other hand, Assumption~\ref{ass:sigma} gives
\[
|\sigma^\omega(t,x)-\sigma^\omega(t,y)|\le K|x-y|^\theta,
\qquad \theta\in[1/2,1],
\]
and
\[
\lambda\le|\sigma^\omega(t,x)|\le\Lambda.
\]
Thus the frozen coefficients satisfy the distribution-free assumptions of
\cite[Theorem~2.3]{HuangWang2023}, and the frozen equation admits a unique strong solution
$X^\omega$.

The uniform linear-growth bound on the drift and the boundedness of the diffusion coefficient
imply, by the Burkholder--Davis--Gundy and Young inequalities, that
\[
\mathbb E\left[\sup_{0\le r\le t}|X_r^\omega|^p\right]
\le
C\bigl(1+\mathbb E[|\xi|^p]\bigr)
+
C\int_0^t
\mathbb E\left[\sup_{0\le \ell\le s}|X_\ell^\omega|^p\right]ds.
\]
Gronwall's lemma yields the moment estimate in \eqref{eq:3.7}, with a constant independent
of $\omega$.

\textup{(ii)} The coefficient bounds in
Assumptions~\ref{ass:b} and~\ref{ass:sigma} are uniform in the frozen
path $\omega$. Thus the transition-density bounds recorded in
\eqref{eq:frozen-gaussian-two-sided}--\eqref{eq:frozen-second-derivative-kernel}
apply uniformly in $\omega$. In particular,
$p^\omega(s,x;t,y)$ is continuous and strictly positive, and
\[
P_{s,t}^\omega f(x)
=\int_{\mathbb R}f(y)p^\omega(s,x;t,y)\,dy.
\]
Integrating \eqref{eq:frozen-gradient-kernel} with respect to $y$ yields
\eqref{eq:3.8}. The asserted second-derivative estimate follows directly from
\eqref{eq:frozen-second-derivative-kernel}, whose constants are uniform in $\omega$.

\textup{(iii)} For $t\in(0,T]$,
\[
u_t^\omega(y)=\int_{\mathbb R}p^\omega(0,x;t,y)u_0(x)\,dx.
\]
The continuity and positivity of $p^\omega$, together with $u_0>0$, imply that
$u_t^\omega$ is continuous and strictly positive.

We next prove the density estimate \eqref{eq:3.10}. Since
\[
\|u_t^\omega\|_{C_b^\gamma}
=\|u_t^\omega\|_\infty+[u_t^\omega]_\gamma,
\]
it suffices to estimate these two terms separately.

We first estimate $\|u_t^\omega\|_\infty$. Using the Gaussian upper bound \eqref{eq:frozen-gaussian-two-sided}
and the uniform bi-Lipschitz estimate \eqref{eq:flow-bilipschitz},
we obtain
\[
\begin{aligned}
u_t^\omega(y)
&=\int_{\mathbb R}p^\omega(0,x;t,y)u_0(x)\,dx\\
&\le C\|u_0\|_\infty
\int_{\mathbb R}g_c\bigl(t,\chi_{t,0}^\omega(x)-y\bigr)\,dx\\
&\le C\|u_0\|_\infty\int_{\mathbb R}g_c(t,z)\,dz\\
&\le C.
\end{aligned}
\]
Therefore,
\begin{equation}
\sup_{\omega\in C_T}\sup_{0\le t\le T}
\|u_t^\omega\|_\infty\le C.
\label{eq:uniform-linfty}
\end{equation}

We next estimate the H\"older seminorm, defined by
\[
[u_t^\omega]_\gamma
=\sup_{y\ne y'}
\frac{|u_t^\omega(y)-u_t^\omega(y')|}{|y-y'|^\gamma}.
\]
It is therefore enough to prove
\[
|u_t^\omega(y)-u_t^\omega(y')|\le C|y-y'|^\gamma.
\]
Following the forward Duhamel representation and parametrix expansion in
\cite[Section~2.1, (2.13), (2.17), and Section~2.2.1, (2.18)--(2.21)]{MenozziPesceZhang2021}, we first work with spatially mollified
coefficients. The mollification index is omitted below. All estimates in this argument
are uniform with respect to this index. Set
\[
a^\omega(r,z):=\bigl(\sigma^\omega(r,z)\bigr)^2,
\]
and define
\[
\begin{aligned}
L_r^\omega f(z)
&=b^\omega(r,z)f'(z)+\frac12a^\omega(r,z)f''(z),\\
\bar L_{r,t}^{\omega,y}f(z)
&=b^\omega\bigl(r,\chi_{r,t}^\omega(y)\bigr)f'(z)
+\frac12a^\omega\bigl(r,\chi_{r,t}^\omega(y)\bigr)f''(z).
\end{aligned}
\]
Set
\[
A_{s,t}^{\omega,y}
:=\int_s^t a^\omega\bigl(r,\chi_{r,t}^\omega(y)\bigr)\,dr,
\qquad
B_{s,t}^{\omega,y}
:=\int_s^t b^\omega\bigl(r,\chi_{r,t}^\omega(y)\bigr)\,dr.
\]
The associated Gaussian density is
\[
\bar p^\omega(s,x;t,y)
=\frac{1}{\sqrt{2\pi A_{s,t}^{\omega,y}}}
\exp\!\left(
-\frac{\bigl(y-x-B_{s,t}^{\omega,y}\bigr)^2}
{2A_{s,t}^{\omega,y}}
\right).
\]
Define
\[
H^\omega(r,z;t,y)
:=\bigl(L_r^\omega-\bar L_{r,t}^{\omega,y}\bigr)
\bar p^\omega(r,z;t,y).
\]
Then the Duhamel--parametrix representation gives
\[
\begin{aligned}
p^\omega(s,x;t,y)
&=
\bar p^\omega(s,x;t,y)\\
&\quad+
\int_s^t\int_{\mathbb R}
p^\omega(s,x;r,z)H^\omega(r,z;t,y)\,dz\,dr.
\end{aligned}
\]

For the frozen coefficients, the parametrix H\"older argument uses the following bounds for the initial Gaussian term and the parametrix kernel. To obtain the first bound, we write the term as a Gaussian mollification of $u_0$. For $|y-y'|\le1$, the uniform bi-Lipschitz estimate \eqref{eq:flow-bilipschitz} and the spatial regularity of the coefficients give increments of its mean and variance bounded by $Ct|y-y'|^\varrho$. Assumption~\ref{ass:init} and $\gamma<\varrho$ then yield a H\"older bound that is uniform for $0<t\le T$, while the case $|y-y'|>1$ follows from the uniform $L^\infty$-bound. For the second bound, we estimate separately the increment coming from the coefficient centers and that coming from the Gaussian kernel. The first uses the spatial $\varrho$-H\"older regularity, and the second uses interpolation between the standard pointwise and Lipschitz parametrix estimates. Splitting according to $|y-y'|\le\sqrt{t-r}$ and $|y-y'|>\sqrt{t-r}$ gives the displayed estimate. The exponent $\varrho-\gamma$ accounts for the spatial regularity of both the drift and the diffusion.
\begin{equation}
\left|
\int_{\mathbb R}
\bigl[\bar p^\omega(0,x;t,y)-\bar p^\omega(0,x;t,y')\bigr]u_0(x)\,dx
\right|
\le C|y-y'|^\gamma.
\label{eq:parametrix-main-holder}
\end{equation}
\begin{equation}
\begin{split}
&|H^\omega(r,z;t,y)-H^\omega(r,z;t,y')|\\
&\quad\le C|y-y'|^\gamma(t-r)^{-1+(\varrho-\gamma)/2}\\
&\qquad\times
\Bigl[
 g_c\bigl(t-r,\chi_{r,t}^\omega(y)-z\bigr)
 +g_c\bigl(t-r,\chi_{r,t}^\omega(y')-z\bigr)
\Bigr].
\end{split}
\label{eq:parametrix-kernel-holder}
\end{equation}

Integrating the parametrix representation against the initial density $u_0$ and applying
Fubini's theorem, we obtain
\[
\begin{aligned}
u_t^\omega(y)
&=\int_{\mathbb R}u_0(x)p^\omega(0,x;t,y)\,dx\\
&=\int_{\mathbb R}u_0(x)\bar p^\omega(0,x;t,y)\,dx
 +\int_0^t\int_{\mathbb R}u_r^\omega(z)H^\omega(r,z;t,y)\,dz\,dr.
\end{aligned}
\]
Since $u_r^\omega$ is a probability density with a uniform $L^\infty$-bound by \eqref{eq:uniform-linfty}, and the integral of $g_c(t-r,\cdot)$ is independent of its spatial center, for $0\le r<t$ and $y,y'\in\mathbb R$ we have
\[
\begin{aligned}
\int_{\mathbb R}u_r^\omega(z)
&\Bigl[
g_c\bigl(t-r,\chi_{r,t}^\omega(y)-z\bigr)
+g_c\bigl(t-r,\chi_{r,t}^\omega(y')-z\bigr)
\Bigr]\,dz\\
&\le \|u_r^\omega\|_\infty
\int_{\mathbb R}\Bigl[
g_c\bigl(t-r,\chi_{r,t}^\omega(y)-z\bigr)
+g_c\bigl(t-r,\chi_{r,t}^\omega(y')-z\bigr)
\Bigr]\,dz\\
&\le C.
\end{aligned}
\]
Hence, combining this estimate with the parametrix estimates
\eqref{eq:parametrix-main-holder} and \eqref{eq:parametrix-kernel-holder}, for
$y,y'\in\mathbb R$, we obtain
\[
\begin{aligned}
|u_t^\omega(y)-u_t^\omega(y')|
&\le
\left|
\int_{\mathbb R}
[\bar p^\omega(0,x;t,y)-\bar p^\omega(0,x;t,y')]u_0(x)\,dx
\right|\\
&\quad+
\int_0^t\int_{\mathbb R}
u_r^\omega(z)
|H^\omega(r,z;t,y)-H^\omega(r,z;t,y')|\,dz\,dr\\
&\le C|y-y'|^\gamma\\
&\quad+C|y-y'|^\gamma
\int_0^t(t-r)^{-1+(\varrho-\gamma)/2}
\int_{\mathbb R}u_r^\omega(z)\\
&\qquad\times\Bigl[
 g_c\bigl(t-r,\chi_{r,t}^\omega(y)-z\bigr)
 +g_c\bigl(t-r,\chi_{r,t}^\omega(y')-z\bigr)
 \Bigr]\,dz\,dr\\
&\le C|y-y'|^\gamma
\left[1+\int_0^t(t-r)^{-1+(\varrho-\gamma)/2}\,dr\right]\\
&=C|y-y'|^\gamma
\left[1+\frac{2}{\varrho-\gamma}t^{(\varrho-\gamma)/2}\right]\\
&\le C|y-y'|^\gamma.
\end{aligned}
\]
Here the last two inequalities use $0<\gamma<\varrho$ and $t\le T$.
Since the preceding estimates are uniform with respect to the mollification
parameter, we may pass to the original coefficients by the approximation
procedure described above. Thus,
\[
\sup_{\omega\in C_T}\sup_{0<t\le T}
[u_t^\omega]_\gamma\le C.
\]
Assumption~\ref{ass:init} supplies the estimate at $t=0$. Combining this H\"older bound
with the estimate \eqref{eq:uniform-linfty}, we obtain
\[
\sup_{\omega\in C_T}\sup_{0\le t\le T}
\|u_t^\omega\|_{C_b^\gamma}\le C,
\]
which proves the density estimate \eqref{eq:3.10}.

It remains to prove the $L^1$-continuity of the frozen density flow. Fix $s\in[0,T]$ and
let $t_n\to s$. Since $X^\omega$ has continuous sample paths,
\[
X_{t_n}^\omega\longrightarrow X_s^\omega
\qquad\text{a.s.},
\]
and hence
\[
\mathcal L(X_{t_n}^\omega)\Longrightarrow\mathcal L(X_s^\omega).
\]
By the uniform $C_b^\gamma$-estimate \eqref{eq:3.10}, the family
$\{u_t^\omega:t\in[0,T]\}$ is uniformly bounded and equicontinuous. Hence, by the
Arzel\`a--Ascoli theorem and a diagonal argument, every subsequence of $(u_{t_n}^\omega)$
admits a further subsequence converging locally uniformly on $\mathbb R$ to a continuous
limit. The local uniform convergence together with the weak convergence above identifies
this limit with $u_s^\omega$.

By the moment estimate \eqref{eq:3.7} and Markov's inequality,
\[
\sup_{0\le t\le T}\int_{|x|>R}u_t^\omega(x)\,dx
\le \frac{C}{R^p},
\qquad R>0.
\]
Along this further subsequence, still denoted by $(t_n)$,
\[
\begin{aligned}
\|u_{t_n}^\omega-u_s^\omega\|_{L^1}
&\le\int_{|x|\le R}|u_{t_n}^\omega(x)-u_s^\omega(x)|\,dx\\
&\quad+\int_{|x|>R}\bigl(u_{t_n}^\omega(x)+u_s^\omega(x)\bigr)\,dx.
\end{aligned}
\]
Consequently,
\[
\limsup_{n\to\infty}\|u_{t_n}^\omega-u_s^\omega\|_{L^1}
\le \frac{2C}{R^p}.
\]
Letting $R\to\infty$, we conclude that this further subsequence converges to $u_s^\omega$
in $L^1$. Since every subsequence admits a further subsequence converging to $u_s^\omega$
in $L^1$, the whole sequence converges. Hence,
\[
u^\omega\in C([0,T];L^1(\mathbb R)).
\]
\end{proof}

\begin{remark}
The uniformity in the frozen path is essential for the contraction argument below and follows
from the uniformity with respect to the second variable in Assumptions~\ref{ass:b} and
\ref{ass:sigma}.
\end{remark}

\begin{proposition}[Uniform Density-Class Inclusion]\label{prop:density-class}
Suppose that Assumptions~\ref{ass:b}--\ref{ass:init} hold. Then there exist constants
\[
K_0>0,\qquad
\delta_0>0,\qquad
0<\varepsilon_0<\min\{\alpha,1-\alpha\}
\]
that are independent of $t\in[0,T]$ and $\omega\in C_T$. Moreover,
\[
u_t^\omega\in
\mathcal S_{\alpha,K_0,\delta_0,\varepsilon_0},
\qquad
t\in[0,T],\quad \omega\in C_T.
\]
\end{proposition}

\begin{proof}
Fix
\[
0<\varepsilon_0<\min\{\alpha,1-\alpha\}.
\]
By the moment estimate \eqref{eq:3.7} in
Proposition~\ref{prop:frozen-density}\textup{(i)} and Markov's inequality,
\[
\int_{|x|\ge K_0}u_t^\omega(x)\,dx
=
\mathbb P(|X_t^\omega|\ge K_0)
\le
\frac{\mathbb E[|X_t^\omega|^p]}{K_0^p}
\le
\frac{C(1+\mathbb E[|\xi|^p])}{K_0^p}.
\]
Hence, choosing $K_0$ sufficiently large, we obtain
\begin{equation}
\int_{|x|\ge K_0}u_t^\omega(x)\,dx\le\varepsilon_0,
\label{eq:3.12}
\end{equation}
uniformly in $t\in[0,T]$ and $\omega\in C_T$. It follows that
\[
F_{u_t^\omega}(-K_0)\le\varepsilon_0<\alpha,
\qquad
F_{u_t^\omega}(K_0)\ge1-\varepsilon_0>\alpha,
\]
and hence
\begin{equation}
Q_\alpha(u_t^\omega)\in[-K_0,K_0].
\label{eq:3.13}
\end{equation}

It remains to prove a uniform positive lower bound on this interval. For the lower heat-kernel estimate, we now choose $\chi_{r,s}^\omega$ as a fixed Peano solution of \eqref{eq:peano-flow}; uniqueness of the flow is not required. By \cite[Remark~1.3]{MenozziPesceZhang2021}, the lower bound remains valid with this center and with constants independent of $\omega$. The estimate is
\begin{equation}
p^\omega(s,x;t,y)
\ge
\frac{c}{\sqrt{t-s}}
\exp\!\left(-\frac{C|\chi_{t,s}^\omega(x)-y|^2}{t-s}\right),
\qquad 0\le s<t\le T.
\label{eq:3.14}
\end{equation}
Fix $t\in(0,T]$ and $|y|\le K_0$. Set
\[
R_0:=K_0+\frac{\sqrt T}{2}.
\]
If $|x-y|\le\sqrt t/2$, then $|x|\le R_0$. Since $u_0$ is continuous and strictly
positive,
\[
m_0:=\inf_{|x|\le R_0}u_0(x)>0.
\]
Since
\[
|b(x,z)|\le C(1+|x|),
\qquad x,z\in\mathbb R,
\]
Gronwall's lemma yields
\[
\sup_{\omega\in C_T}\sup_{0\le r\le T}\sup_{|x|\le R_0}
|\chi_{r,0}^\omega(x)|\le C.
\]
Consequently,
\[
|\chi_{t,0}^\omega(x)-x|
\le
\int_0^t|b(\chi_{r,0}^\omega(x),\omega_r)|\,dr
\le Ct.
\]
For $|x-y|\le\sqrt t/2$,
\[
|\chi_{t,0}^\omega(x)-y|
\le Ct+\frac{\sqrt t}{2},
\]
and therefore
\[
\frac{|\chi_{t,0}^\omega(x)-y|^2}{t}
\le
2C^2t+\frac12
\le
2C^2T+\frac12.
\]
The lower heat-kernel estimate in \eqref{eq:3.14} therefore yields a constant $c_0>0$,
uniform over
\[
t\in(0,T],\qquad |y|\le K_0,\qquad |x-y|\le\sqrt t/2,\qquad \omega\in C_T,
\]
such that
\[
p^\omega(0,x;t,y)\ge\frac{c_0}{\sqrt t}.
\]
Using the density representation and restricting the integral to
$|x-y|\le\sqrt t/2$, we obtain
\[
\begin{aligned}
u_t^\omega(y)
&\ge
\int_{|x-y|\le\sqrt t/2}p^\omega(0,x;t,y)u_0(x)\,dx\\
&\ge
m_0\frac{c_0}{\sqrt t}
\int_{|x-y|\le\sqrt t/2}dx\\
&=m_0c_0.
\end{aligned}
\]
At $t=0$, the continuity and strict positivity of $u_0$ give
\[
\inf_{|y|\le K_0}u_0(y)>0.
\]
Thus, with
\[
\delta_0
:=
\min\left\{
\inf_{|y|\le K_0}u_0(y),\,m_0c_0
\right\}>0,
\]
we obtain
\begin{equation}
\inf_{\omega\in C_T}\inf_{0\le t\le T}\inf_{|y|\le K_0}
u_t^\omega(y)
\ge\delta_0.
\label{eq:uniform-density-lower-bound}
\end{equation}
Combining \eqref{eq:3.12}, \eqref{eq:3.13}, and
\eqref{eq:uniform-density-lower-bound}, we obtain
\[
u_t^\omega\in
\mathcal S_{\alpha,K_0,\delta_0,\varepsilon_0},
\qquad
t\in[0,T],\quad \omega\in C_T.
\]
This completes the proof.
\end{proof}

\subsection{Stability of Frozen Densities}

Let $\zeta$ be an $\mathcal F_0$-measurable random variable such that
\[
\mathbb E[|\zeta|^p]<\infty
\]
for some $p\ge2$, and suppose that its law admits a density
\[
u_0^\zeta\in\mathcal D\cap C_b^\gamma(\mathbb R).
\]
Let $\omega\in C_T$. We denote by $X^{\omega,\zeta}$ the solution of
\[
 dX_t^{\omega,\zeta}
 =b(X_t^{\omega,\zeta},\omega_t)\,dt
 +\sigma(X_t^{\omega,\zeta},\omega_t)\,dW_t,
 \qquad X_0^{\omega,\zeta}=\zeta,
\]
and by $u_t^{\omega,\zeta}$ its density. When the initial condition is the fixed random
variable $\xi$ from Assumption~\ref{ass:init}, we continue to write $X^\omega$ and
$u_t^\omega$.

The conclusions of Proposition~\ref{prop:frozen-density} remain valid for
$X^{\omega,\zeta}$, with the same proof.
In particular, the moment and $C_b^\gamma$-bounds used below have constants depending on
$\mathbb E[|\zeta|^p]$ and $\|u_0^\zeta\|_{C_b^\gamma}$, but independent of $\omega$.

By Assumption~\ref{ass:b},
\begin{equation}
|b(x,\omega_t)-b(x,\widetilde\omega_t)|
\le
C|\omega_t-\widetilde\omega_t|,
\qquad x\in\mathbb R,
\label{eq:3.15}
\end{equation}
where $C>0$ is independent of $x$, $t$, $\omega$, and $\widetilde\omega$.

Set
\[
a(x,z):=\sigma^2(x,z).
\]
Assumption~\ref{ass:sigma} and $\gamma<\theta$ give
\[
\|a(\cdot,z)-a(\cdot,\widetilde z)\|_\infty
\le
C|z-\widetilde z|.
\]
The same assumptions and the spatial $\theta$-H\"older continuity also give
\[
[\sigma(\cdot,z)+\sigma(\cdot,\widetilde z)]_\gamma\le C.
\]
Using the mixed increment condition, we obtain
\[
\begin{aligned}
[a(\cdot,z)-a(\cdot,\widetilde z)]_\gamma
&\le
2\Lambda[\sigma(\cdot,z)-\sigma(\cdot,\widetilde z)]_\gamma\\
&\quad+
\|\sigma(\cdot,z)-\sigma(\cdot,\widetilde z)\|_\infty
[\sigma(\cdot,z)+\sigma(\cdot,\widetilde z)]_\gamma\\
&\le
C|z-\widetilde z|.
\end{aligned}
\]
Hence
\begin{equation}
\|a(\cdot,\omega_t)-a(\cdot,\widetilde\omega_t)\|_{C_b^\gamma}
\le
C|\omega_t-\widetilde\omega_t|.
\label{eq:3.16}
\end{equation}

For the spatial translation estimate below, set
\[
\bar\gamma:=\frac{p\gamma}{p+1}\in(0,\gamma),
\]
where $p\ge2$ denotes the moment exponent associated with the initial law.

The drift contribution in the frozen-density comparison can be controlled directly by the gradient bound \eqref{eq:3.8}. The diffusion contribution instead involves $\partial_{xx}P_{s,t}^\omega f$, whose kernel bound has the non-integrable time singularity $(t-s)^{-1}$. The following lemma uses the zero-integral property of $\partial_{xx}p^\omega$ and the H\"older regularity of the coefficient-density product to reduce this singularity to the integrable order $(t-s)^{-1+\bar\gamma/2}$. This is precisely the estimate needed for the diffusion term in Proposition~\ref{prop:frozen-stability}.

\begin{lemma}[Second-Order Perturbation Estimate]
\label{lem:second-order-perturbation}
Suppose that Assumptions~\ref{ass:b} and~\ref{ass:sigma} hold.
Let $p\ge2$, let $\eta$ be an $\mathcal F_0$-measurable random
variable such that $\mathbb E[|\eta|^p]<\infty$, and suppose that
its law admits a density $u_0^\eta\in\mathcal D\cap C_b^\gamma(\mathbb R)$.
Then, for every
$\omega,\widetilde\omega\in C_T$ and $0\le s<t\le T$,
\begin{equation}
\begin{aligned}
&\left\|\int_{\mathbb R}
\bigl[a(x,\omega_s)-a(x,\widetilde\omega_s)\bigr]
 u_s^{\widetilde\omega,\eta}(x)
 \partial_{xx}p^\omega(s,x;t,\cdot)\,dx
 \right\|_{L^1}\\
&\qquad\le C(t-s)^{-1+\bar\gamma/2}
 |\omega_s-\widetilde\omega_s|.
\end{aligned}
\label{eq:second-order-kernel}
\end{equation}
The constant $C$ is independent of $\omega,\widetilde\omega,s,t$
and may depend on $T$, the coefficient bounds, the $p$-moment of $\eta$,
and $\|u_0^\eta\|_{C_b^\gamma}$.
\end{lemma}

\begin{proof}
We divide the proof into four steps and first work with spatially smooth
approximations. All bounds below are uniform in the approximation parameter.

\noindent\textbf{Step 1: Estimate the coefficient--density product.} We derive H\"older and weighted $L^1$-bounds, which will control its spatial translations. By \eqref{eq:3.10}, \eqref{eq:3.16}, and the moment estimate \eqref{eq:3.7}
applied with initial variable $\eta$, we have
\begin{align}
&\Big[\bigl(a(\cdot,\omega_s)-a(\cdot,\widetilde\omega_s)\bigr)
u_s^{\widetilde\omega,\eta}\Big]_\gamma\notag\\
&\le
\|a(\cdot,\omega_s)-a(\cdot,\widetilde\omega_s)\|_\infty
[u_s^{\widetilde\omega,\eta}]_\gamma\notag\\
&\quad+\|u_s^{\widetilde\omega,\eta}\|_\infty
[a(\cdot,\omega_s)-a(\cdot,\widetilde\omega_s)]_\gamma\notag\\
&\le C|\omega_s-\widetilde\omega_s|,\label{eq:holder-product}\\
\intertext{and}
&\int_{\mathbb R}(1+|x|^p)
|a(x,\omega_s)-a(x,\widetilde\omega_s)|
u_s^{\widetilde\omega,\eta}(x)\,dx\notag\\
&\le C|\omega_s-\widetilde\omega_s|
\int_{\mathbb R}(1+|x|^p)u_s^{\widetilde\omega,\eta}(x)\,dx\notag\\
&\le C|\omega_s-\widetilde\omega_s|.\label{eq:weighted-l1}
\end{align}

\medskip
\begingroup\allowdisplaybreaks[1]
\noindent\textbf{Step 2: Control spatial translations.} We convert the H\"older and weighted bounds into an $L^1$-translation estimate with a smaller H\"older exponent. Define
\[
H_s(x):=[a(x,\omega_s)-a(x,\widetilde\omega_s)]
u_s^{\widetilde\omega,\eta}(x).
\]
For $0<|w|\leq 1$ and $M\geq 2$, \eqref{eq:holder-product}--\eqref{eq:weighted-l1} give
\begin{align*}
\int_{\mathbb R}|H_s(x+w)-H_s(x)|\,dx
&\leq 2M[H_s]_\gamma|w|^\gamma
+\int_{|x|>M}|H_s(x+w)|\,dx
+\int_{|x|>M}|H_s(x)|\,dx\\
&\leq C|\omega_s-\widetilde\omega_s|
\bigl(M|w|^\gamma+M^{-p}\bigr).
\end{align*}
For the tail bound, the change of variables $y=x+w$ and $|w|\leq1$ show that both integration regions are contained in $\{|z|>M/2\}$, since $M\geq2$. The estimate \eqref{eq:weighted-l1} then provides the factor $M^{-p}$. Choosing
\[
M=2|w|^{-\gamma/(p+1)}
\]
and using $\bar\gamma=p\gamma/(p+1)$ yield
$M|w|^\gamma=2|w|^{\bar\gamma}$ and $M^{-p}=2^{-p}|w|^{\bar\gamma}$. For $|w|>1$, \eqref{eq:weighted-l1} gives the same bound directly, and the case $w=0$ is trivial. Consequently,
\begin{equation}
\int_{\mathbb R}|H_s(x+w)-H_s(x)|\,dx
\leq C|\omega_s-\widetilde\omega_s||w|^{\bar\gamma},
\qquad w\in\mathbb R.
\label{eq:3.19}
\end{equation}
\endgroup

\medskip
\noindent\textbf{Step 3: Center the perturbation at the Gaussian point.} We prove the zero-integral identity for the second derivative of the kernel, which allows us to subtract the product at the kernel center. By the bi-Lipschitz estimate \eqref{eq:flow-bilipschitz}, we have
\[
|\chi_{t,s}^\omega(x)-y|
\geq C^{-1}
\left|x-(\chi_{t,s}^\omega)^{-1}(y)\right|.
\]
Combining this with the second-derivative Gaussian
bound \eqref{eq:frozen-second-derivative-kernel}, possibly with a different $c>0$,
gives
\[
\begin{aligned}
|\partial_{xx}p^\omega(s,x;t,y)|
\leq C(t-s)^{-1}
g_c\left(t-s,
x-(\chi_{t,s}^\omega)^{-1}(y)\right).
\end{aligned}
\]
Consequently,
\[
\begin{aligned}
\int_{\mathbb R}
|\partial_{xx}p^\omega(s,x;t,y)|\,dx
&\leq C(t-s)^{-1}
\int_{\mathbb R}
g_c\left(t-s,
x-(\chi_{t,s}^\omega)^{-1}(y)\right)\,dx\\
&\leq C(t-s)^{-1}<\infty.
\end{aligned}
\]

On the other hand, by the first-derivative Gaussian
bound \eqref{eq:frozen-gradient-kernel} and the bi-Lipschitz estimate \eqref{eq:flow-bilipschitz},
\[
\begin{aligned}
|\partial_xp^\omega(s,x;t,y)|
&\leq C(t-s)^{-1/2}
g_c(t-s,\chi_{t,s}^\omega(x)-y)\\
&=C(t-s)^{-1}
\exp\left(
-\frac{|\chi_{t,s}^\omega(x)-y|^2}
{c(t-s)}
\right)
\longrightarrow 0
\end{aligned}
\]
as $|x|\to\infty$.

Therefore, by the fundamental theorem of calculus,
for every $R>0$,
\[
\int_{-R}^{R}
\partial_{xx}p^\omega(s,x;t,y)\,dx
=
\partial_xp^\omega(s,R;t,y)
-\partial_xp^\omega(s,-R;t,y).
\]
Together with the integrability established above,
letting $R\to\infty$ yields
\[
\begin{aligned}
\int_{\mathbb R}
\partial_{xx}p^\omega(s,x;t,y)\,dx
&=\lim_{R\to\infty}
\int_{-R}^{R}
\partial_{xx}p^\omega(s,x;t,y)\,dx\\
&=\lim_{R\to\infty}
\left[
\partial_xp^\omega(s,R;t,y)
-\partial_xp^\omega(s,-R;t,y)
\right]\\
&=0.
\end{aligned}
\]

Since the Gaussian bound for $\partial_{xx}p^\omega(s,x;t,y)$ is centered at
$x=(\chi_{t,s}^\omega)^{-1}(y)$, we subtract the coefficient-density product
evaluated at this point. Using the zero-integral identity, for every
$y\in\mathbb R$, we obtain
\begin{equation}
\begin{aligned}
&\int_{\mathbb R}
[a(x,\omega_s)-a(x,\widetilde\omega_s)]
u_s^{\widetilde\omega,\eta}(x)
\partial_{xx}p^\omega(s,x;t,y)\,dx\\
&=\int_{\mathbb R}
\Big\{
[a(x,\omega_s)-a(x,\widetilde\omega_s)]
u_s^{\widetilde\omega,\eta}(x)\\
&\qquad-
[a((\chi_{t,s}^\omega)^{-1}(y),\omega_s)
-a((\chi_{t,s}^\omega)^{-1}(y),\widetilde\omega_s)]
u_s^{\widetilde\omega,\eta}
((\chi_{t,s}^\omega)^{-1}(y))
\Big\}\\
&\qquad\qquad\times
\partial_{xx}p^\omega(s,x;t,y)\,dx.
\end{aligned}
\label{eq:cancellation-identity}
\end{equation}

\medskip
\noindent\textbf{Step 4: Obtain an integrable time singularity.} We combine the translation estimate with Gaussian integration to complete the perturbation bound. The spatially regularized drift has a bounded
continuous spatial derivative. Hence its flow is
spatially differentiable, and \eqref{eq:flow-bilipschitz} yields
\[
C^{-1}
\leq |\partial_z\chi_{t,s}^\omega(z)|
\leq C,
\qquad z\in\mathbb R,\quad 0\leq s\leq t\leq T.
\]

Combining \eqref{eq:cancellation-identity} with the Gaussian upper bound
for the second spatial derivative of the transition density in
\eqref{eq:frozen-second-derivative-kernel}, we first change variables
$y=\chi_{t,s}^\omega(z)$ using the uniform Jacobian bounds and then set $w=x-z$.
The bi-Lipschitz estimate also allows us to replace the Gaussian center
$\chi_{t,s}^\omega(x)-\chi_{t,s}^\omega(z)$ by $x-z$ after enlarging the constant
$c$ in the Gaussian bound. Applying the spatial translation estimate
\eqref{eq:3.19} to the resulting inner integral, we obtain
\[
\begin{aligned}
&\left\|
\int_{\mathbb R}
[a(x,\omega_s)-a(x,\widetilde\omega_s)]
u_s^{\widetilde\omega,\eta}(x)
\partial_{xx}p^\omega(s,x;t,\cdot)\,dx
\right\|_{L^1}\\
&\leq C(t-s)^{-1}
\int_{\mathbb R}\int_{\mathbb R}
\Big|
[a(x,\omega_s)-a(x,\widetilde\omega_s)]
u_s^{\widetilde\omega,\eta}(x)\\
&\qquad-
[a(z,\omega_s)-a(z,\widetilde\omega_s)]
u_s^{\widetilde\omega,\eta}(z)
\Big|
g_c(t-s,x-z)\,dx\,dz\\
&=C(t-s)^{-1}
\int_{\mathbb R}g_c(t-s,w)
\int_{\mathbb R}
\Big|
[a(z+w,\omega_s)-a(z+w,\widetilde\omega_s)]
u_s^{\widetilde\omega,\eta}(z+w)\\
&\qquad-
[a(z,\omega_s)-a(z,\widetilde\omega_s)]
u_s^{\widetilde\omega,\eta}(z)
\Big|\,dz\,dw\\
&\leq C|\omega_s-\widetilde\omega_s|
(t-s)^{-1}
\int_{\mathbb R}|w|^{\bar\gamma}
g_c(t-s,w)\,dw\\
&=C|\omega_s-\widetilde\omega_s|
(t-s)^{-1+\bar\gamma/2}
\int_{\mathbb R}|z|^{\bar\gamma}g_c(1,z)\,dz\\
&\leq C|\omega_s-\widetilde\omega_s|
(t-s)^{-1+\bar\gamma/2}.
\end{aligned}
\]

Since $\bar\gamma>0$, the time singularity $(t-s)^{-1+\bar\gamma/2}$ is integrable as $s\uparrow t$. The approximation argument above then yields the estimate for the original coefficients.
\end{proof}

\begin{proposition}[Stability of Frozen Densities]\label{prop:frozen-stability}
Suppose that Assumptions~\ref{ass:b} and~\ref{ass:sigma} hold. Let $\xi$ and $\eta$ be
$\mathcal F_0$-measurable
random variables defined on the same probability space. Suppose that
\[
\mathbb E[|\xi|^p]+\mathbb E[|\eta|^p]<\infty
\]
for some $p\ge2$, and that their laws admit densities
\[
u_0^\xi,u_0^\eta\in\mathcal D\cap C_b^\gamma(\mathbb R).
\]
Then there exists a constant $C>0$ such that, for every
$\omega,\widetilde\omega\in C_T$ and every $t\in(0,T]$,
\begin{equation}
\begin{aligned}
\|u_t^{\omega,\xi}-u_t^{\widetilde\omega,\eta}\|_{L^1}
\le{}&
Ct^{-1/2}\mathbb E[|\xi-\eta|]\\
&+C\int_0^t(t-s)^{-1/2}|\omega_s-\widetilde\omega_s|\,ds\\
&+C\int_0^t(t-s)^{-1+\bar\gamma/2}|\omega_s-\widetilde\omega_s|\,ds.
\end{aligned}
\label{eq:3.17}
\end{equation}
If the diffusion coefficient is independent of the second variable, the last integral is absent.
The constant $C$ is independent of $\omega$ and $\widetilde\omega$. It depends only on $T$
and the constants in Assumptions~\ref{ass:b} and~\ref{ass:sigma}. It may also depend on the
relevant moment bounds and the $C_b^\gamma$-bounds of the initial densities.
\end{proposition}

\begin{proof}
Fix $t\in(0,T]$, and let
$f\in C_c^\infty(\mathbb R)$ satisfy
$\|f\|_{L^\infty}\leq 1$.
We first establish the duality identity for spatially
smooth approximations of the coefficients.
For notational simplicity, the approximation index
is omitted throughout the following computation.

For $0\leq s\leq t$, define the operator
$L_s^\omega$ and its formal adjoint by
\[
L_s^\omega\varphi(x)
:=
b(x,\omega_s)\varphi'(x)
+\frac12a(x,\omega_s)\varphi''(x),
\]
and
\[
(L_s^\omega)^*\psi(x)
:=
-\partial_x\big(b(x,\omega_s)\psi(x)\big)
+\frac12\partial_{xx}
\big(a(x,\omega_s)\psi(x)\big).
\]

Set
\[
v(s,x):=P_{s,t}^\omega f(x).
\]

Then $v$ solves the backward Kolmogorov equation
\[
\left\{
\begin{aligned}
\partial_s v(s,x)+L_s^\omega v(s,x)&=0,
&&0\leq s<t,\quad x\in\mathbb R,\\
v(t,x)&=f(x),
&&x\in\mathbb R.
\end{aligned}
\right.
\]

The corresponding density $u_s^{\omega,\xi}$
solves the forward Kolmogorov equation
\[
\left\{
\begin{aligned}
\partial_s u_s^{\omega,\xi}(x)
&=(L_s^\omega)^*u_s^{\omega,\xi}(x),
&&0<s\leq t,\quad x\in\mathbb R,\\
u_0^{\omega,\xi}(x)&=u_0^\xi(x),
&&x\in\mathbb R.
\end{aligned}
\right.
\]

By the backward and forward equations and the
adjoint relation,
\[
\begin{aligned}
\frac{d}{ds}
\int_{\mathbb R}
v(s,x)u_s^{\omega,\xi}(x)\,dx
&=
\int_{\mathbb R}
\partial_s v(s,x)u_s^{\omega,\xi}(x)\,dx\\
&\quad+
\int_{\mathbb R}
v(s,x)\partial_s u_s^{\omega,\xi}(x)\,dx\\
&=
-\int_{\mathbb R}
L_s^\omega v(s,x)u_s^{\omega,\xi}(x)\,dx\\
&\quad+
\int_{\mathbb R}
v(s,x)(L_s^\omega)^*
u_s^{\omega,\xi}(x)\,dx\\
&=0.
\end{aligned}
\]

Integrating over $[0,t]$ and using
$v(t,x)=f(x)$, we obtain
\[
\begin{aligned}
\int_{\mathbb R}
f(x)u_t^{\omega,\xi}(x)\,dx
&=
\int_{\mathbb R}
v(0,x)u_0^\xi(x)\,dx\\
&=\mathbb E[v(0,\xi)].
\end{aligned}
\]

Since $u_s^{\widetilde\omega,\eta}$
satisfies the forward equation associated with
$L_s^{\widetilde\omega}$, we have
\[
\begin{aligned}
\frac{d}{ds}
\int_{\mathbb R}
v(s,x)u_s^{\widetilde\omega,\eta}(x)\,dx
&=
-\int_{\mathbb R}
L_s^\omega v(s,x)
u_s^{\widetilde\omega,\eta}(x)\,dx\\
&\quad+
\int_{\mathbb R}
v(s,x)(L_s^{\widetilde\omega})^*
u_s^{\widetilde\omega,\eta}(x)\,dx\\
&=
\int_{\mathbb R}
(L_s^{\widetilde\omega}-L_s^\omega)
v(s,x)
u_s^{\widetilde\omega,\eta}(x)\,dx.
\end{aligned}
\]

Integrating over $[0,t]$ and using
$v(t,x)=f(x)$, we obtain
\[
\begin{aligned}
\int_{\mathbb R}
f(x)u_t^{\widetilde\omega,\eta}(x)\,dx
&=
\mathbb E[v(0,\eta)]\\
&\quad+
\int_0^t\int_{\mathbb R}
(L_s^{\widetilde\omega}-L_s^\omega)
v(s,x)
u_s^{\widetilde\omega,\eta}(x)
\,dx\,ds.
\end{aligned}
\]

Subtracting this identity from the corresponding
identity for $u_t^{\omega,\xi}$, we obtain
\[
\begin{aligned}
&\int_{\mathbb R}
f(x)\big(
u_t^{\omega,\xi}(x)
-u_t^{\widetilde\omega,\eta}(x)
\big)\,dx\\
&=
\mathbb E[v(0,\xi)-v(0,\eta)]\\
&\quad+
\int_0^t\int_{\mathbb R}
(L_s^\omega-L_s^{\widetilde\omega})
v(s,x)
u_s^{\widetilde\omega,\eta}(x)
\,dx\,ds.
\end{aligned}
\]

By the definition of the generators,
\[
\begin{aligned}
&(L_s^\omega-L_s^{\widetilde\omega})
v(s,x)\\
&=
[b(x,\omega_s)-b(x,\widetilde\omega_s)]
\partial_xv(s,x)\\
&\quad+
\frac12
[a(x,\omega_s)-a(x,\widetilde\omega_s)]
\partial_{xx}v(s,x).
\end{aligned}
\]

Since $v(s,x)=P_{s,t}^\omega f(x)$, we have
\[
\partial_xv(s,x)
=\partial_xP_{s,t}^\omega f(x),
\qquad
\partial_{xx}v(s,x)
=\partial_{xx}P_{s,t}^\omega f(x).
\]

Therefore,
\begin{equation}
\begin{aligned}
&\int_{\mathbb R}
f(x)\big(
u_t^{\omega,\xi}(x)
-u_t^{\widetilde\omega,\eta}(x)
\big)\,dx\\
&=I_0+I_1+I_2,
\end{aligned}
\label{eq:3.18}
\end{equation}
where
\[
I_0:=\mathbb E[v(0,\xi)-v(0,\eta)],
\]
\[
\begin{aligned}
I_1:=
\int_0^t\int_{\mathbb R}
&[b(x,\omega_s)-b(x,\widetilde\omega_s)]
\partial_xP_{s,t}^\omega f(x)\\
&\times
u_s^{\widetilde\omega,\eta}(x)
\,dx\,ds,
\end{aligned}
\]
and
\[
\begin{aligned}
I_2:=
\frac12\int_0^t\int_{\mathbb R}
&[a(x,\omega_s)-a(x,\widetilde\omega_s)]
\partial_{xx}P_{s,t}^\omega f(x)\\
&\times
u_s^{\widetilde\omega,\eta}(x)
\,dx\,ds.
\end{aligned}
\]

We next estimate the three terms separately.
By the mean value theorem and the gradient estimate
\eqref{eq:3.8}, we have
\[
\begin{aligned}
|I_0|
&=
|\mathbb E[v(0,\xi)-v(0,\eta)]|\\
&\leq
\|\partial_xv(0,\cdot)\|_\infty
\mathbb E[|\xi-\eta|]\\
&\leq
Ct^{-1/2}\mathbb E[|\xi-\eta|].
\end{aligned}
\]

Next, by the drift estimate \eqref{eq:3.15} and the gradient
estimate \eqref{eq:3.8},
\[
\begin{aligned}
|I_1|
&\leq
C\int_0^t
(t-s)^{-1/2}|\omega_s-\widetilde\omega_s|
\left(
\int_{\mathbb R}
u_s^{\widetilde\omega,\eta}(x)\,dx
\right)ds\\
&=
C\int_0^t
(t-s)^{-1/2}
|\omega_s-\widetilde\omega_s|\,ds.
\end{aligned}
\]

For $I_2$, the kernel representation
\[
\partial_{xx}P_{s,t}^\omega f(x)
=
\int_{\mathbb R}
\partial_{xx}p^\omega(s,x;t,y)f(y)\,dy
\]
and Fubini's theorem, together with Lemma~\ref{lem:second-order-perturbation},
yield
\[
\begin{aligned}
|I_2|
&\leq
\frac12\|f\|_\infty
\int_0^t
\left\|
\int_{\mathbb R}
[a(x,\omega_s)-a(x,\widetilde\omega_s)]
u_s^{\widetilde\omega,\eta}(x)
\partial_{xx}p^\omega(s,x;t,\cdot)\,dx
\right\|_{L^1}
ds\\
&\leq
C\|f\|_\infty
\int_0^t
(t-s)^{-1+\bar\gamma/2}
|\omega_s-\widetilde\omega_s|\,ds\\
&\leq
C\int_0^t
(t-s)^{-1+\bar\gamma/2}
|\omega_s-\widetilde\omega_s|\,ds.
\end{aligned}
\]

The estimates above are uniform in the approximation
parameter. Passing to the limit by the spatial
approximation procedure described above yields
the same bounds for the original coefficients.

Combining the estimates for $I_0$, $I_1$, and $I_2$,
we obtain
\[
\begin{aligned}
&\left|
\int_{\mathbb R}
f(x)\big(
u_t^{\omega,\xi}(x)-u_t^{\widetilde\omega,\eta}(x)
\big)\,dx
\right|\\
&\leq
Ct^{-1/2}\mathbb E[|\xi-\eta|]\\
&\quad+
C\int_0^t
(t-s)^{-1/2}
|\omega_s-\widetilde\omega_s|\,ds\\
&\quad+
C\int_0^t
(t-s)^{-1+\bar\gamma/2}
|\omega_s-\widetilde\omega_s|\,ds.
\end{aligned}
\]

For $h\in L^1(\mathbb R)$, the $L^1$-duality formula
gives
\[
\|h\|_{L^1}
=
\sup_{\substack{f\in C_c^\infty(\mathbb R)\\
\|f\|_{L^\infty}\leq1}}
\left|
\int_{\mathbb R}f(x)h(x)\,dx
\right|.
\]

Taking
\[
h=u_t^{\omega,\xi}-u_t^{\widetilde\omega,\eta}
\]
yields the stability estimate in \eqref{eq:3.17}.

If the diffusion coefficient is independent of the
second variable, then
\[
a(x,\omega_s)-a(x,\widetilde\omega_s)=0,
\]
so $I_2=0$ and the last integral in \eqref{eq:3.17} is absent.
\end{proof}

\subsection{Quantile Stability and the Frozen Quantile Map}

We next prove the stability of the quantile map on the density class introduced in
Section~\ref{sec:model}.

\begin{lemma}[Quantile Stability]\label{lem:quantile-stability}
Let $u,v\in\mathcal S_{\alpha,K_0,\delta_0,\varepsilon_0}$. Then
\begin{equation}
|Q_\alpha(u)-Q_\alpha(v)|
\le
\delta_0^{-1}\|u-v\|_{L^1}.
\label{eq:3.20}
\end{equation}
\end{lemma}

\begin{proof}
Write
\[
q_u=Q_\alpha(u),
\qquad
q_v=Q_\alpha(v).
\]
Both quantiles belong to $[-K_0,K_0]$. Suppose, without loss of generality, that
$q_u\le q_v$. Since $v\ge\delta_0$ on $[-K_0,K_0]$,
\[
\delta_0(q_v-q_u)
\le
\int_{q_u}^{q_v}v(x)\,dx.
\]
Moreover, $F_u(q_u)=F_v(q_v)=\alpha$, and hence
\[
\begin{aligned}
\int_{q_u}^{q_v}v(x)\,dx
&=F_v(q_v)-F_v(q_u)\\
&=F_u(q_u)-F_v(q_u)\\
&=\int_{-\infty}^{q_u}(u(x)-v(x))\,dx\\
&\le\|u-v\|_{L^1}.
\end{aligned}
\]
Therefore, $q_v-q_u\le\delta_0^{-1}\|u-v\|_{L^1}$. Interchanging $u$ and $v$ gives the
same estimate when $q_v\le q_u$, proving the quantile stability estimate
\eqref{eq:3.20}.
\end{proof}

Recall from Section~\ref{sec:model} that, for $\omega\in C_T$, the frozen quantile path is
given by
\[
(\Phi(\omega))_t=Q_\alpha(u_t^\omega),\qquad t\in[0,T].
\]
By Propositions~\ref{prop:frozen-density} and \ref{prop:density-class}, together with Lemma~\ref{lem:quantile-stability}, for every $s,t\in[0,T]$,
\[
|Q_\alpha(u_t^\omega)-Q_\alpha(u_s^\omega)|
\le
\delta_0^{-1}\|u_t^\omega-u_s^\omega\|_{L^1}.
\]
Hence, $t\mapsto Q_\alpha(u_t^\omega)$ is continuous, and therefore
$\Phi:C_T\to C_T$ is well defined.

We next apply Proposition~\ref{prop:frozen-stability} to two frozen equations with the same
initial condition. Combining the resulting estimate with
Lemma~\ref{lem:quantile-stability} yields stability of the frozen quantile map.

\begin{proposition}[Stability of the Frozen Quantile Map]\label{prop:phi}
There exists a constant $C_\Phi>0$ such that, for every
$\omega,\widetilde\omega\in C_T$ and every $t\in(0,T]$,
\begin{equation}
\begin{aligned}
|(\Phi(\omega))_t-(\Phi(\widetilde\omega))_t|
\le{}&
C_\Phi\int_0^t(t-s)^{-1/2}|\omega_s-\widetilde\omega_s|\,ds\\
&+C_\Phi\int_0^t(t-s)^{-1+\bar\gamma/2}
|\omega_s-\widetilde\omega_s|\,ds.
\end{aligned}
\label{eq:3.21}
\end{equation}
The constant $C_\Phi$ is independent of $\omega$ and $\widetilde\omega$.
\end{proposition}

\begin{proof}
At $t=0$, both frozen densities equal $u_0$, so the quantile difference is zero and the estimate is trivial. For $t>0$, applying Proposition~\ref{prop:frozen-stability} with the same initial condition $\xi$ gives
\[
\begin{aligned}
\|u_t^\omega-u_t^{\widetilde\omega}\|_{L^1}
\le{}&
C\int_0^t(t-s)^{-1/2}|\omega_s-\widetilde\omega_s|\,ds\\
&+C\int_0^t(t-s)^{-1+\bar\gamma/2}
|\omega_s-\widetilde\omega_s|\,ds.
\end{aligned}
\]
Consequently, Lemma~\ref{lem:quantile-stability} yields
\[
\begin{aligned}
|(\Phi(\omega))_t-(\Phi(\widetilde\omega))_t|
&=|Q_\alpha(u_t^\omega)-Q_\alpha(u_t^{\widetilde\omega})|\\
&\le\delta_0^{-1}\|u_t^\omega-u_t^{\widetilde\omega}\|_{L^1}\\
&\le
C_\Phi\int_0^t(t-s)^{-1/2}|\omega_s-\widetilde\omega_s|\,ds\\
&\quad+
C_\Phi\int_0^t(t-s)^{-1+\bar\gamma/2}
|\omega_s-\widetilde\omega_s|\,ds.
\end{aligned}
\]
This proves the frozen-quantile stability estimate \eqref{eq:3.21}.
\end{proof}

\begin{remark}
Taking the supremum in the stability estimate \eqref{eq:3.21} gives
\[
\|\Phi(\omega)-\Phi(\widetilde\omega)\|_{C([0,t])}
\le
C_\Phi\left(2\sqrt t+\frac{2}{\bar\gamma}t^{\bar\gamma/2}\right)
\|\omega-\widetilde\omega\|_{C([0,t])},
\qquad t\in[0,T].
\]
Hence, on a sufficiently short time interval $[0,t]$, the frozen quantile map is a contraction
under the usual supremum norm. This recovers the local fixed-point mechanism used by
Kolokoltsov~\cite{Kolokoltsov2013}. In the proof of Theorem~\ref{thm:wellposed}, we instead
use an equivalent exponentially
weighted supremum norm to obtain a contraction directly on the whole interval $[0,T]$.
\end{remark}

\subsection{Proof of Theorem~\ref{thm:wellposed}}

For $\kappa>0$, define the exponentially weighted supremum norm by
\[
\|\omega\|_{\kappa,T}
:=
\sup_{0\le t\le T}e^{-\kappa t}|\omega_t|,
\qquad \omega\in C_T.
\]
Since
\[
e^{-\kappa T}\|\omega\|_T
\le
\|\omega\|_{\kappa,T}
\le
\|\omega\|_T,
\]
this norm is equivalent to the usual supremum norm. Hence,
$(C_T,\|\cdot\|_{\kappa,T})$ is a Banach space.

By Proposition~\ref{prop:phi}, for every $t\in(0,T]$,
\[
\begin{aligned}
&e^{-\kappa t}|(\Phi(\omega))_t-(\Phi(\widetilde\omega))_t|\\
&\quad\le
C_\Phi\|\omega-\widetilde\omega\|_{\kappa,T}
\left[
\int_0^t e^{-\kappa(t-s)}(t-s)^{-1/2}\,ds
+
\int_0^t e^{-\kappa(t-s)}(t-s)^{-1+\bar\gamma/2}\,ds
\right]\\
&\quad\le
C_\Phi\left[
\sqrt\pi\,\kappa^{-1/2}
+\Gamma(\bar\gamma/2)\kappa^{-\bar\gamma/2}
\right]
\|\omega-\widetilde\omega\|_{\kappa,T}.
\end{aligned}
\]
Here we used
\[
\int_0^\infty e^{-\kappa r}r^{-1/2}\,dr
=
\sqrt\pi\,\kappa^{-1/2},
\qquad
\int_0^\infty e^{-\kappa r}r^{-1+\bar\gamma/2}\,dr
=
\Gamma(\bar\gamma/2)\kappa^{-\bar\gamma/2}.
\]
Therefore,
\begin{equation}
\|\Phi(\omega)-\Phi(\widetilde\omega)\|_{\kappa,T}
\le
C_\Phi\left[
\sqrt\pi\,\kappa^{-1/2}
+\Gamma(\bar\gamma/2)\kappa^{-\bar\gamma/2}
\right]
\|\omega-\widetilde\omega\|_{\kappa,T}.
\label{eq:3.23}
\end{equation}
Choose $\kappa>0$ sufficiently large so that
\[
C_\Phi\left[
\sqrt\pi\,\kappa^{-1/2}
+\Gamma(\bar\gamma/2)\kappa^{-\bar\gamma/2}
\right]<1.
\]
Then the map $\Phi$ is a strict contraction. Banach's fixed-point theorem therefore yields a
unique $q\in C_T$ such that
\begin{equation}
q=\Phi(q).
\label{eq:3.24}
\end{equation}
Equivalently, $q_t=Q_\alpha(u_t^q)$ for every $t\in[0,T]$.

Let $X^q$ be the unique strong solution of the frozen equation
\[
 dX_t^q=b(X_t^q,q_t)\,dt+\sigma(X_t^q,q_t)\,dW_t,
 \qquad X_0^q=\xi.
\]
By Proposition~\ref{prop:frozen-density}, every $X_t^q$ has a density $u_t^q\in\mathcal D$. Since $q_t=Q_\alpha(u_t^q)$ by \eqref{eq:3.24}, $X^q$ satisfies Eq.~\eqref{eq:2.1}. Its quantile path is the continuous path $q$. Thus $X:=X^q$ is a strong solution in the class specified in Section~\ref{sec:model}.

Let $X$ and $\widetilde X$ be two strong solutions that have the same initial condition and
are driven by the same Brownian motion. Denote their density and quantile paths by $(u,q)$
and $(\widetilde u,\widetilde q)$, respectively. Let $X^q$ and $X^{\widetilde q}$ denote the
unique strong solutions of the frozen equations corresponding to $q$ and $\widetilde q$,
respectively. Since each pair solves the same frozen equation with the same initial condition
and Brownian motion, pathwise uniqueness gives, up to indistinguishability,
\[
X=X^q,\qquad \widetilde X=X^{\widetilde q}.
\]
Hence
\[
q=\Phi(q),\qquad \widetilde q=\Phi(\widetilde q).
\]
By uniqueness of the fixed point, $q=\widetilde q$. Pathwise uniqueness for the corresponding
frozen equation then yields that $X$ and $\widetilde X$ are indistinguishable.

Finally, the uniform estimate in Proposition~\ref{prop:frozen-density} gives
\[
\mathbb E\left[\sup_{0\le t\le T}|X_t|^p\right]
=
\mathbb E\left[\sup_{0\le t\le T}|X_t^q|^p\right]
\le
C\bigl(1+\mathbb E[|\xi|^p]\bigr).
\]
This proves the moment estimate \eqref{eq:2.8} and completes the proof.
\qed

\subsection{Proof of Theorem~\ref{thm:stability}}

Let $X^\xi$ and $X^\eta$ be the solutions with initial data $\xi$ and $\eta$, driven by the
same Brownian motion. Denote their densities by $u_t^\xi$ and $u_t^\eta$, and their quantile
paths by
\[
q_t^\xi=Q_\alpha(u_t^\xi),
\qquad
q_t^\eta=Q_\alpha(u_t^\eta).
\]
Set
\[
Z_t=X_t^\xi-X_t^\eta,
\qquad
R(t)=\|u_t^\xi-u_t^\eta\|_{L^1}.
\]
By pathwise uniqueness for the frozen equations, $X^\xi$ and $X^{q^\xi,\xi}$ are
indistinguishable, and so are $X^\eta$ and $X^{q^\eta,\eta}$. Hence
\[
u_t^\xi=u_t^{q^\xi,\xi},
\qquad
u_t^\eta=u_t^{q^\eta,\eta},
\qquad t\in[0,T].
\]

\noindent\textbf{Step 1: Derive a Volterra inequality for the density difference.} We combine frozen-density stability with quantile stability so that the quantile difference no longer appears explicitly on the right-hand side. Because $\sigma(x,z)=\sigma(x)$ in Theorem~\ref{thm:stability}, the diffusion perturbation
vanishes and the last integral in Proposition~\ref{prop:frozen-stability} is absent. Applying
Proposition~\ref{prop:frozen-stability} with
\[
\omega=q^\xi,
\qquad
\widetilde\omega=q^\eta
\]
gives, for every $t\in(0,T]$,
\begin{equation}
R(t)
\le
Ct^{-1/2}\mathbb E[|\xi-\eta|]
+C\int_0^t\frac{|q_s^\xi-q_s^\eta|}{\sqrt{t-s}}\,ds.
\label{eq:3.25}
\end{equation}
Fix
\[
\varepsilon_0\in(0,\min\{\alpha,1-\alpha\}).
\]
Applying the tail estimate of Proposition~\ref{prop:density-class} to the two initial laws,
choose $K_0>0$ such that, uniformly in $t\in[0,T]$,
\[
\int_{|x|\ge K_0}u_t^\xi(x)\,dx\le\varepsilon_0,
\qquad
\int_{|x|\ge K_0}u_t^\eta(x)\,dx\le\varepsilon_0.
\]
Consequently,
\[
|q_t^\xi|\le K_0,
\qquad
|q_t^\eta|\le K_0,
\qquad t\in[0,T].
\]
With the interval $[-K_0,K_0]$ fixed, we apply the positive lower-bound argument in
Proposition~\ref{prop:density-class} separately to $u_0^\xi$ and $u_0^\eta$. This gives
constants $\delta_\xi>0$ and $\delta_\eta>0$ such that
\[
\inf_{0\le t\le T}\inf_{|x|\le K_0}u_t^\xi(x)\ge\delta_\xi,
\qquad
\inf_{0\le t\le T}\inf_{|x|\le K_0}u_t^\eta(x)\ge\delta_\eta.
\]
Set
\[
\delta_0:=\min\{\delta_\xi,\delta_\eta\}>0.
\]
Then, for every $t\in[0,T]$,
\[
u_t^\xi,u_t^\eta\in\mathcal S_{\alpha,K_0,\delta_0,\varepsilon_0}.
\]
Therefore, Lemma~\ref{lem:quantile-stability} gives
\begin{equation}
|q_t^\xi-q_t^\eta|
\le
\delta_0^{-1}\|u_t^\xi-u_t^\eta\|_{L^1}
=
\delta_0^{-1}R(t).
\label{eq:3.26}
\end{equation}
Combining the density stability estimate \eqref{eq:3.25} and the quantile stability estimate
\eqref{eq:3.26}, we obtain
\begin{equation}
R(t)
\le
Ct^{-1/2}\mathbb E[|\xi-\eta|]
+C\int_0^t\frac{R(s)}{\sqrt{t-s}}\,ds.
\label{eq:3.27}
\end{equation}

\noindent\textbf{Step 2: Iterate the Volterra inequality.} We use the beta-function identity to sum the iterated forcing terms and show that the remainder vanishes, yielding the density and quantile bounds. Define, for $t>0$,
\[
F(t)=t^{-1/2}\mathbb E[|\xi-\eta|],
\qquad
(Jh)(t)=\int_0^t\frac{h(s)}{\sqrt{t-s}}\,ds.
\]
Fix $C_V>0$ such that the Volterra inequality \eqref{eq:3.27} can be written as
\[
R(t)\le C_VF(t)+C_V(JR)(t).
\]
Iterating this inequality gives, for every $m\ge0$,
\[
R(t)
\le
C_V\sum_{k=0}^m C_V^k(J^kF)(t)
+C_V^{m+1}(J^{m+1}R)(t).
\]
The beta-function identity, valid for $a>0$, gives
\[
\begin{aligned}
\int_0^t(t-s)^{-1/2}s^{a-1}\,ds
&=t^{a-1/2}B\!\left(\frac12,a\right)\\
&=t^{a-1/2}
\frac{\Gamma(1/2)\Gamma(a)}{\Gamma(a+1/2)}.
\end{aligned}
\]
Repeated iteration gives, for every $k\ge0$,
\[
(J^kF)(t)
=
\frac{\pi^{(k+1)/2}}{\Gamma((k+1)/2)}
 t^{k/2-1/2}\mathbb E[|\xi-\eta|].
\]
Therefore, for $0<t\le T$,
\[
C_V^{k+1}(J^kF)(t)
\le
Ct^{-1/2}
\frac{(C_V\sqrt{\pi T})^k\sqrt\pi}{\Gamma((k+1)/2)}
\mathbb E[|\xi-\eta|].
\]
Since
\[
\sum_{k=0}^\infty
\frac{(C_V\sqrt{\pi T})^k}{\Gamma((k+1)/2)}<\infty,
\]
there exists $C>0$ such that
\[
C_V\sum_{k=0}^\infty C_V^k(J^kF)(t)
\le
Ct^{-1/2}\mathbb E[|\xi-\eta|].
\]
Moreover, since $u_t^\xi$ and $u_t^\eta$ are probability densities, $R(t)\le2$. The
iterated-kernel formula therefore gives
\[
C_V^{m+1}(J^{m+1}R)(t)
\le
2\frac{(C_V\sqrt\pi)^{m+1}t^{(m+1)/2}}
{\Gamma((m+1)/2+1)}
\longrightarrow0
\]
as $m\to\infty$, uniformly for $t\in[0,T]$. Letting $m\to\infty$, we obtain
\begin{equation}
R(t)
\le
Ct^{-1/2}\mathbb E[|\xi-\eta|],
\qquad t\in(0,T].
\label{eq:3.28}
\end{equation}
That is,
\[
\|u_t^\xi-u_t^\eta\|_{L^1}
\le
Ct^{-1/2}\mathbb E[|\xi-\eta|].
\]
The quantile stability estimate \eqref{eq:3.26} now gives
\begin{equation}
|q_t^\xi-q_t^\eta|
\le
Ct^{-1/2}\mathbb E[|\xi-\eta|],
\qquad t\in(0,T].
\label{eq:3.29}
\end{equation}

\noindent\textbf{Step 3: Estimate the difference of the solution processes.} We apply the Yamada--Watanabe approximation and Gronwall's lemma, using the quantile bound obtained in Step 2. Recall that $Z_t=X_t^\xi-X_t^\eta$. Then
\[
 dZ_t
 =\bigl[b(X_t^\xi,q_t^\xi)-b(X_t^\eta,q_t^\eta)\bigr]dt
 +\bigl[\sigma(X_t^\xi)-\sigma(X_t^\eta)\bigr]dW_t.
\]
Applying It\^o's formula to $V_{r,\varepsilon}(Z_t)$ and taking expectations gives
\[
\begin{aligned}
\mathbb E[V_{r,\varepsilon}(Z_t)]
&=
\mathbb E[V_{r,\varepsilon}(Z_0)]\\
&\quad+
\mathbb E\int_0^t
V_{r,\varepsilon}'(Z_s)
\bigl[b(X_s^\xi,q_s^\xi)-b(X_s^\eta,q_s^\eta)\bigr]ds\\
&\quad+
\frac12\mathbb E\int_0^t
V_{r,\varepsilon}''(Z_s)
|\sigma(X_s^\xi)-\sigma(X_s^\eta)|^2\,ds.
\end{aligned}
\]
Since
\[
\left|
V_{r,\varepsilon}'(Z_s)
\bigl(\sigma(X_s^\xi)-\sigma(X_s^\eta)\bigr)
\right|
\le 2\Lambda,
\]
the stochastic integral is square-integrable and hence has zero expectation.
Decompose the drift difference as
\[
b(X_s^\xi,q_s^\xi)-b(X_s^\eta,q_s^\eta)
=D_{1,s}+D_{2,s}+D_{3,s},
\]
where
\[
D_{1,s}
=b_1(X_s^\xi,q_s^\xi)-b_1(X_s^\eta,q_s^\xi),
\]
\[
D_{2,s}
=b_1(X_s^\eta,q_s^\xi)-b_1(X_s^\eta,q_s^\eta),
\]
and
\[
D_{3,s}
=b_2(X_s^\xi,q_s^\xi)-b_2(X_s^\eta,q_s^\eta).
\]
By the monotonicity of $b_1$ in the spatial variable and the sign property in
\eqref{eq:3.3},
\[
V_{r,\varepsilon}'(Z_s)D_{1,s}\le0.
\]
By the Lipschitz continuity of $b_1$ in the second variable,
\[
|V_{r,\varepsilon}'(Z_s)D_{2,s}|
\le
C|q_s^\xi-q_s^\eta|.
\]
By the joint Lipschitz continuity of $b_2$,
\[
|V_{r,\varepsilon}'(Z_s)D_{3,s}|
\le
C|Z_s|+C|q_s^\xi-q_s^\eta|.
\]
For the diffusion term, the estimate \eqref{eq:3.5} gives
\[
V_{r,\varepsilon}''(Z_s)
|\sigma(X_s^\xi)-\sigma(X_s^\eta)|^2
\le
\frac{C}{\log r}\varepsilon^{2\theta-1}.
\]
Since $V_{r,\varepsilon}(Z_0)\le|Z_0|=|\xi-\eta|$, we obtain
\[
\begin{aligned}
\mathbb E[V_{r,\varepsilon}(Z_t)]
&\le
\mathbb E[|\xi-\eta|]
+C\int_0^t\mathbb E[|Z_s|]ds
+C\int_0^t|q_s^\xi-q_s^\eta|ds\\
&\quad+
\frac{Ct}{\log r}\varepsilon^{2\theta-1}.
\end{aligned}
\]
The approximation bounds \eqref{eq:3.2} imply
\[
\mathbb E[|Z_t|]
\le
\mathbb E[V_{r,\varepsilon}(Z_t)]+\varepsilon.
\]
If $\theta>1/2$, letting $\varepsilon\downarrow0$ gives
\[
\frac{\varepsilon^{2\theta-1}}{\log r}\longrightarrow0.
\]
If $\theta=1/2$, first letting $\varepsilon\downarrow0$ and then $r\to\infty$ gives
$1/\log r\to0$. Hence, in either case,
\[
\mathbb E[|Z_t|]
\le
\mathbb E[|\xi-\eta|]
+C\int_0^t\mathbb E[|Z_s|]ds
+C\int_0^t|q_s^\xi-q_s^\eta|ds.
\]
Using the quantile stability estimate \eqref{eq:3.29} for $s>0$ and the
integrability of $s^{-1/2}$ at the origin, we obtain
\[
\int_0^t|q_s^\xi-q_s^\eta|ds
\le
C\mathbb E[|\xi-\eta|]\int_0^t s^{-1/2}ds
\le
2C\sqrt T\,\mathbb E[|\xi-\eta|].
\]
Therefore,
\[
\mathbb E[|Z_t|]
\le
C\mathbb E[|\xi-\eta|]
+C\int_0^t\mathbb E[|Z_s|]ds.
\]
Gronwall's lemma yields
\[
\mathbb E[|Z_t|]
\le
C\mathbb E[|\xi-\eta|],
\qquad t\in[0,T].
\]
Thus,
\[
\sup_{0\le t\le T}\mathbb E[|X_t^\xi-X_t^\eta|]
\le
C\mathbb E[|\xi-\eta|].
\]
Combining this estimate with the density stability estimate \eqref{eq:3.28} and the
quantile stability estimate \eqref{eq:3.29}, we obtain the estimates in
\eqref{eq:2.9}--\eqref{eq:2.11}. This completes the proof of
Theorem~\ref{thm:stability}.
\qed

\end{document}